\documentclass[journal]{IEEEtran}
\makeatletter
\def\subsubsection{\@startsection{subsubsection}
                                 {3}
                                 {\z@}
                                 {0ex plus 0.1ex minus 0.1ex}
                                 {0ex}
                                 {\normalfont\normalsize\bfseries}}
\makeatother
\usepackage[T1]{fontenc}

\usepackage{ulem}
\usepackage{amsmath}
\allowdisplaybreaks
\usepackage{hhline}
\usepackage{yfonts,color}
\usepackage{soul,xcolor}
\usepackage{verbatim}
\usepackage{amsmath}
\allowdisplaybreaks
\usepackage{amssymb}
\usepackage{amsthm}
\usepackage{float}
\usepackage{bm}
\usepackage{url}
\usepackage{array}
\usepackage{cite}
\usepackage{graphicx}
\usepackage{framed} 
\usepackage{balance} 
\usepackage{epsfig,epstopdf}
\usepackage{booktabs}
\usepackage{caption}
\usepackage{courier}
\usepackage{pseudocode}

\usepackage{lineno}

\usepackage{enumerate}
\usepackage{algorithm}
\usepackage{algpseudocode}
\newtheorem{definition}{Definition}
\newtheorem{theorem}{Theorem}
\newtheorem{lemma}[theorem]{Lemma}
\newtheorem{proposition}[theorem]{Proposition}
\newtheorem{corollary}[theorem]{Corollary}
\newtheorem{assumption}{Assumption}

\newcommand{\rom}[1]{\uppercase\expandafter{\romannumeral #1\relax}}
\usepackage{color}
\usepackage{soul,xcolor}
\newcommand{\nm}[1]{{\color{blue!20!black!30!green}\bf{[NM: #1]}}}

\usepackage{cancel}

\begin{document} 
\setulcolor{red}
\setul{red}{2pt}
\title{Finite Rate Distributed Weight-Balancing and Average Consensus Over Digraphs}
\author{Chang-Shen Lee, Nicol\`{o} Michelusi, and Gesualdo Scutari

\thanks{Lee and Michelusi are with the School of Electrical and Computer Engineering, Purdue University, West Lafayette, IN, USA; emails:
<lee2495,michelus>@purdue.edu.}
\thanks{Scutari is with the School of Industrial Engineering, Purdue University, West Lafayette, IN, USA; email:
gscutari@purdue.edu. The work of Scutari and part of the work of Lee have  been supported by   the USA National Science Foundation under Grants  CIF 1719205  and CIF 1564044; and in part by the Army Research Office under Grant W911NF1810238.
}
\thanks{A preliminary version of this paper appeared at IEEE CDC 2018 \cite{Lee_2018CDC}.}
\vspace{-0.7cm}
}
\maketitle

\setstcolor{red}

\begin{abstract}
This paper proposes   the first distributed algorithm that solves the weight-balancing problem using only finite rate and simplex communications among nodes, compliant with the directed nature of the graph edges. It is proved that the algorithm converges  to a weight-balanced solution at sublinear rate.
The analysis builds upon a new metric inspired by positional  system representations, which characterizes the dynamics of information exchange over the network, and on a novel step-size rule.
 Building on this result, a novel distributed algorithm is proposed that solves the {\it average} consensus problem over digraphs, using, at each timeslot, finite rate simplex communications between adjacent nodes -- some bits for the weight-balancing problem and  others for the average consensus.   
Convergence of the proposed quantized consensus algorithm to the average of the node's {\it unquantized} initial values is established, both almost surely and in the moment generating function of the error; and  a sublinear convergence rate is proved for sufficiently large step-sizes. Numerical results validate our theoretical findings.
\end{abstract}
\vspace{-2mm}
\begin{IEEEkeywords}
Distributed weight-balancing, distributed average consensus, directed graph, quantization, data rate.
\end{IEEEkeywords}
\vspace{-5mm}
\section{Introduction} \label{section_intro}\vspace{-0.1cm}
Digraphs play a key role in a number of network applications, such as distributed optimization \cite{Lorenzo2016}, distributed flow-balancing \cite{Hadjicostis2019}, distributed averaging and cooperative control \cite{Olfati-Saber2004}, to name a few.  
In particular,  distributed average consensus, whereby  nodes aim  at  agreeing on the sample average of their local values,  has received considerable attention over the  years;   some applications include load-balancing \cite{Cybenko1989}, vehicle formation \cite{Fax2004}, and sensor networks \cite{Scutari08}. 
Several of the these distributed algorithms, when run on digraphs,  
  require some form of graph regularity, such as the \emph{weight-balanced} property \cite{Rikos2014}: at each node, the sum of the outgoing edge weights equals that of the incoming edge weights.

Several centralized algorithms have been proposed to balance a digraph; see, e.g., \cite{Loh1970} and references therein. In this paper, we are interested in  the design of {\it distributed}   algorithms that solve the weight-balancing and average consensus problems over digraphs, using only {\it quantized} information, {\it simplex communications},\footnote{One way, as opposed to duplex, two ways communications.} and without knowledge of the graph topology other than the direct neighbor.
This problem is motivated by realistic scenarios, such as wireless sensor networks, where channels may be asymmetric due to different transmit powers of nodes and interference, and where communications are subject to finite rate constraints.
To date, the design of such algorithms in distributed settings remains a challenging and open problem, as documented next.
 \vspace{-3mm}
\subsection{Related works} \label{Sec:related_work}
\noindent \textbf{Distributed weight-balancing} algorithms were proposed in \cite{Hadjicostis2019,Gharesifard2012,Rikos2014,Rikos2019int,Rikos2018_delay} (see Table \ref{table_wb}). With the exception of \cite[Sec. \rom{4}]{Rikos2014},\cite{Rikos2019int,Rikos2018_delay}, all these algorithms require \emph{infinite bits} in each communication round, since nodes need to exchange either real valued, or integer but unbounded quantities. Although \cite[Sec. \rom{4}]{Rikos2014} and \cite{Rikos2018_delay} use a finite number of bits at each iteration, 
 this number cannot be   arbitrarily chosen (e.g. to satisfy some  transmission constraints), it is instead  the result of the algorithmic trajectory and thus  it is not known a-priori.   
 In addition, these works adopt \emph{unicast} communications, whereby   nodes transmit  different signals to different out-neighbors. 
 To reduce signaling overhead, \emph{broadcast} communications are preferable in dense networks.
Finally, while compliant with prescribed finite rate constraints, the distributed integer weight-balancing algorithm \cite{Rikos2019int} requires {\it full-duplex} edge communications--each node must exchange information with \emph{both its out- and in-neighbors}--which does not comply with simplex constraints. To the best of our knowledge, the design of distributed algorithms that   solve the weight-balancing problem using a prescribed \emph{finite rate and simplex} communications is an open problem.  

\begin{table}
	\small\centering
	\caption{Related works on graph weight-balancing}
	\vspace{-2mm}
	\label{table_wb}
    \begin{tabular}{| >{\centering\arraybackslash}m{0.915in} | >{\centering\arraybackslash}m{0.46in} | >{\centering\arraybackslash}m{0.38in} | >{\centering\arraybackslash}m{1.02in} |}    
    \hline
    \textbf{\bf Reference} & \textbf{Broadcast} & \textbf{Digraph} & \textbf{\# Bits/Timeslot}  \\ \hline
    \cite{Hadjicostis2019} & \checkmark & & Infinite  \\ \hline
    \cite{Gharesifard2012} & & \checkmark & Infinite  \\ \hline
    \cite[Sec. \rom{3}]{Rikos2014} & \checkmark & \checkmark & Infinite  \\ \hline
    \cite{Rikos2019int} & \checkmark &  & Any  \\ \hline
    \cite[Sec. \rom{4}]{Rikos2014},\cite{Rikos2018_delay} & & \checkmark & Problem-Dependent  \\ \hline
    Proposed & \checkmark & \checkmark & Any  \\ \hline
    \end{tabular}\medskip 
    \vspace{-8mm}
\end{table}


\begin{table*}[tb]
	\small\centering
    \begin{tabular}{| >{\centering\arraybackslash}m{0.6in} | >{\centering\arraybackslash}m{0.8in} | >{\centering\arraybackslash}m{0.6in} | >{\centering\arraybackslash}m{0.7in} | >{\centering\arraybackslash}m{0.95in} | >{\centering\arraybackslash}m{1.1in} | >{\centering\arraybackslash}m{0.8in} |}    
    \hline
    \textbf{\bf References} & \textbf{Quantization} & \textbf{Digraph} &  \textbf{Convergence} & \textbf{Limit Point} & \textbf{\# Bits/Timeslot} & \textbf{Initial Value}  \\ \hline
    \cite{Kashyap2007} & Deterministic &  &  & Neighborhood & {Problem-Dependent} & Integer  \\ \hline
    \cite{Nedic2009} & Deterministic &  & \checkmark & Neighborhood & {Problem-Dependent} & Integer  \\ \hline    
    \cite{Lavaei2012} & Deterministic &  &  & Neighborhood & Problem-Dependent & Box  \\ \hline
    \cite{ElChamie2016} & Deterministic &  & & Neighborhood & Infinite & Any  \\ \hline
    \cite{Zhu2016} & Deterministic &  & & Neighborhood & Any & Any \\ \hline
    \cite{Aysal2008} & Probabilistic &  & \checkmark & Neighborhood & Any & Any  \\ \hline
    \cite{Kar2010} & Probabilistic & & \checkmark & Neighborhood & Any & Box \\ \hline
    \cite{Li2011, Thanou2013} & Deterministic &  & \checkmark & Average & Any & Box  \\ \hline
    \cite{Rajagopal2011} & Probabilistic &  & \checkmark & Average & Any & Box  \\ \hline
    \cite{Li2016_tn} & Deterministic & \checkmark (Balanced) & \checkmark & Average & Any & Box  \\ \hline
    \cite{Lee2017conf} & Probabilistic & \checkmark (Balanced) & \checkmark & Average & Any & Informative \\ \hline
    \cite{Wang2011, Chen2019} & Deterministic & \checkmark & \checkmark & Weighted Average & Any & Box \\ \hline
    \cite{Rikos2018, Rikos2019} & Deterministic & \checkmark & \checkmark & Average & Trajectory-Dependent & Integer   \\ \hline
    \cite{Bernadette2018} & Deterministic & \checkmark & & Neighborhood & Infinite & Box  \\ \hline 
    \cite{Zhu2015} & Probabilistic & \checkmark & \checkmark & Average & Infinite & Any  \\ \hline
    Proposed & Probabilistic & \checkmark & \checkmark & Average & Any & Informative \\ \hline
    \end{tabular}
    \caption{Related works on quantized consensus. problem-dependent: depends on the problem setting, e.g., initial values, weight matrix; trajectory-dependent: depends on trajectory of the algorithm during execution; informative: cf. Assumption \ref{assump_sbar}.}\label{table_ac} 
    \vspace{-5mm}    
\end{table*}

\noindent \textbf{Distributed average consensus} algorithms have a long history, tracing back to the seminal works\cite{Tsitsiklis1984,Xiao2004,Olfati-Saber2004}. These early works assumed that nodes can reliably exchange  unquantized information over undirected networks. To cope with limited data rates, quantization was later introduced, and its effect analyzed for both undirected \cite{Kashyap2007,Nedic2009,Lavaei2012,ElChamie2016,Zhu2016,Aysal2008,Kar2010,Rajagopal2011,Li2011,Thanou2013} 
or directed graphs \cite{Lee2017conf,Wang2011,Li2016_tn, Chen2019, Rikos2018, Rikos2019, Bernadette2018, Zhu2015}, as documented in Table \ref{table_ac}. The quantized average consensus problem based on deterministic uniform quantization and dithered (probabilistic) quantization has been  considered in \cite{Kashyap2007, Nedic2009, Lavaei2012, ElChamie2016, Zhu2016}   and  \cite{Aysal2008,Kar2010}, respectively. However, these schemes do not achieve exact consensus but converge to a neighborhood of the average. \emph{Exact} average consensus is proved in
   \cite{Li2011,Thanou2013} for deterministic quantization and in \cite{Rajagopal2011}
 for probabilistic quantization. However, all these algorithms
  consider \emph{undirected} graphs, which can be easily weight-balanced (e.g. using the Metropolis weights \cite{Xiao2005}). While the extensions of the above deterministic and probabilistic schemes to digraphs were studied in \cite{Wang2011, Li2016_tn, Chen2019} and \cite{Lee2017conf}, respectively, all these works only achieve exact average convergence over \emph{balanced} digraph. However, the weight matrices of digraphs are inherently unbalanced, thus requiring specific weight-balancing algorithms, as documented earlier; they thus  suffer from the same limitations of distributed weight balancing schemes.\vspace{-0.1cm} 

To address unbalanced digraphs, the idea adopted in the seminal work \cite{Kempe2003} is to estimate and compensate the bias caused by the unbalanced weights, via the so-called \emph{push-sum algorithm} \cite{Kempe2003}. This algorithm requires unquantized communication. Unfortunately, applying naively a  finite-bit quantization to the push-sum scheme does not lead to convergence, as we will demonstrate numerically in Section \ref{section_simulation} (cf. Q-Push-Sum). Extensions of push-sum employing   quantization have been developed in \cite{Rikos2018, Rikos2019, Bernadette2018, Zhu2015}. However, \cite{Bernadette2018,Zhu2015} consider  unbounded quantization intervals, which necessitates {\it infinite} bits to encode the signal whereas  \cite{Rikos2018, Rikos2019} impose integer constraints on the initial values of the consensus signals  and necessitate a trajectory-dependent  number of bits. Besides quantization, other instances of imperfect communications in average consensus problems over digraphs were  investigated in \cite{Blasa2011,Olshevsky2018} (asynchrony) and \cite{Fagnani2009,Gerencser2019} (link failure).


\vspace{-3mm}

\subsection{Summary of the main contributions} \label{Sec:contribution}
The above  literature review shows that there are no  distributed algorithms  solving the weight-balancing and the {\it exact} average consensus problems for real initial values over digraphs, using {\it finite-bit quantized information with a prescribed number of bits} and simplex communications. 
This paper provides an answer to these open questions. 

\noindent\textbf{1) Distributed  quantized weight-balancing:} The first contribution is a novel distributed  quantized weight-balancing algorithm whereby nodes transfer part of their balance--the difference between the out-going and the incoming sum-weights, which should be zero for a weight-balanced graph--to their out-neighbors via quantized simplex communications; by doing so, the balance is transferred from high imbalance to low imbalance nodes, provably converging to a weight-balanced solution at sublinear rate. Differently from existing quantized weight-balancing schemes\cite{Gharesifard2012,Rikos2014,Rikos2018_delay}, the proposed algorithm can use at each iteration  a prescribed number of bits (possibly, time-varying). The convergence analysis is also a novel technical contribution of the paper:

\begin{enumerate} 
\item[i)] First,
we identify necessary and sufficient conditions under which the total imbalance decreases, denoted by the \emph{decreasing event} (see D.\ref{def_event_D}). Roughly speaking, this event occurs when a node transfers its balance to a neighbor with balance of opposite sign. Hence, nodes closer to nodes with balance of opposite sign more directly contribute to trigger the {\it decreasing event} and thus reduce  the total imbalance, and are therefore more \emph{important} than those farther away. 
\item[ii)] The next step is to prove that the decreasing events occur often enough that the total imbalance asymptotically vanishes at sublinear rate. To this end, we show that the time interval between two consecutive occurrences of a decreasing event is  {\it uniformly bounded}. This is proved by introducing   a sophisticated metric, a non-negative integer-valued function of the imbalances of nodes and of their importance, which {\it strictly} increases every time there is a transfer of balance from less important nodes to more important ones, {\it up until}  the next decreasing event occurs. By proving that this function is uniformly bounded, we conclude that the decreasing events occur infinitely often.

 To build such a function, we  use the idea of positional system representation:     the   value of the function at each timeslot is expressed by a number whose $h$th digit represents the sum-imbalance of the $h$th most important \emph{nodes}. By doing so,   every transfer of balance from   nodes of lower importance towards those of higher importance causes this function to increase, as it induces a ``carry'' operation from a digit to the next more significant one in its positional representation. 
\item[iii)]  We introduce a novel diminishing step-size rule, which guarantees that the balance at each node is expressed as an integer multiple of the current step-size. 
 This choice greatly facilitates the convergence analysis, since it allows one to tightly control the amount of decrement of the total imbalance at each timeslot.
\end{enumerate}

\noindent\textbf{2) Distributed average quantized consensus:} Building on the proposed weight-balancing scheme, we   introduce a novel distributed algorithm that performs average consensus and weight-balancing 
\emph{on the same time scale} with finite-bit simplex communications--some bits for consensus and some to balance the digraph. For instance, one may perform  one-bit (simplex) communication per channel use, by exchanging weight-balancing and consensus information alternately.
 The key idea behind the algorithm is to preserve the average of the  variables over time, while gradually weight-balancing the graph. 
We prove convergence of the nodes' local variables to the \emph{exact} average of the initial values, both almost surely and in the moment generating function of the error. A sublinear convergence rate is proved for sufficiently large step-sizes.

The rest of the paper is organized as follows. In Section~\ref{sec:background}, we introduce some basic notation and preliminary definitions. Section \ref{section_summary} introduces the ideas of the proposed distributed quantized weight-balancing and average consensus algorithms, whose details are discussed in Section \ref{section_weight_balance} and Section~\ref{section_consensus}, respectively. Some numerical results are discussed in Section~\ref{section_simulation}, while  Section~\ref{section_conclusion} draws some conclusions. The proof of auxiliary lemmas is provided in the appendix. 
\vspace{-3mm}
\section{Notation and Background}\label{sec:background}\vspace{-2mm}
\subsection{Notation}  \label{sec:notation}
The sets of real, integer, nonnegative integer, and positive integer numbers are denoted by $\mathbb{R}$, $\mathbb{Z}$, $\mathbb{Z}_+$, and $\mathbb{Z}_{++}$, respectively. The indicator function is denoted by $\mathcal{I}\{\mathcal{A}\}$,with $\mathcal{I}\{\mathcal{A}\}{=}1$ if $\mathcal{A}$ is true, and $\mathcal{I}\{\mathcal{A}\}{=}0$ otherwise. We define the floor and ceiling functions $\lfloor x\rfloor{\triangleq} \max\{y{\in}\mathbb Z{:} y{\leq}x\}$, $\lceil x\rceil{\triangleq}\min\{y{\in}\mathbb Z{:}y{\geq}x\}$, the sign of $x$, ${\rm sgn}(x){=}x/|x|,\forall x{\neq}0$, ${\rm sgn}(0){=}0$, and the clip function ${\rm clip}_{[l,u]}(x){=}\min\left\{\max\left\{x,l\right\},u\right\}$. We adopt the big-O notation, $f(x){=}\mathcal O\left(g(x)\right){\Leftrightarrow} \limsup_{x \to \infty} \frac{\vert f(x)\vert}{g(x)}{<}\infty$.
Vectors are denoted as $\mathbf x$ (lowercase, boldface), matrices as $\mathbf X$ (uppercase, boldface). All equalities and inequalities involving random variables are tacitly assumed to hold almost surely (i.e., with probability $1$), unless otherwise stated. We use L.x, C.x, D.x, T.x, P.x, A.x and App.x for Lemma x, Corollary x, Definition x, Theorem x, Proposition x, Assumption x and Appendix x, respectively.
 The rest of the symbols used in the paper are summarized in Table~\ref{tableI}.\vspace{-0.2cm}
	\begin{table}
	\small\centering
    \begin{tabular}{| >{\centering\arraybackslash}m{1.15in} | >{\centering\arraybackslash}m{2in} |}    
    \hline
    {\bf Symbol} & {\bf Description} \\ \hline
    $\mathcal{N}_i^+,\mathcal{N}_i^-$ & Out- ($+$) \& in-neighbors ($-$) of node $i$  \\ \hline
    $d_i^+,d_i^-$ & Out- ($+$) \& in-degrees ($-$) of node $i$  \\ \hline
    $k$ & Timeslot Index \\ \hline
    $\gamma(k)$ & Step-size (for weight-balancing)  \\ \hline
    $\alpha(k)$ & Step-size (for consensus)  \\ \hline
    ${\bf A}(k) = \left(a_{ij}(k)\right)_{i,j=1}^N$ & Weight matrix  \\ \hline    
    $S_i^+(k),S_i^-(k)$ & Sum of outgoing ($+$) \& incoming ($-$) weights at node $i$  \\ \hline
    ${\bf b}(k) = \left(b_i(k)\right)_{i=1}^N$ & (Weight) balance  \\ \hline
    ${\bf L}^+(k),{\bf L}^-(k)$ & Graph Laplacian matrices  \\ \hline
    $[q_{\min},q_{\max}]$ & Quantization range for consensus \\ \hline
    ${\bf y}(k){=}\left(y_i(k)\right)_{i=1}^N$ & Local estimate (cf. (\ref{update_consensus}))  \\ \hline
    $\tilde{\bf y}(k){=}\left(\tilde{y}_i(k)\right)_{i=1}^N$ & Clipped local estimate (cf. (\ref{clippedY})) \\ \hline  
    ${\bf y}(0) = \left(y_i(0)\right)_{i=1}^N$ & Initial measurements \\ \hline
    $B_i^{(w)}(k)$ & Number of bits to quantize $b_i(k)$ \\ \hline   
    $B_i^{(c)}(k)$ & Number of bits to quantize $y_i(k)$ \\ \hline  
    $\Delta(B)$ & Distance between consecutive quantization points \\ \hline  
    $n_i(k)$ & Weight-balancing signal sent by node $i$ \\ \hline  
    $\mathcal D_k, \mathcal U_k$ & Decreasing \& Update events \\ \hline 
    \end{tabular}\medskip \caption{Notation used in the paper}\label{tableI}\vspace{-0.3cm}
	\end{table}

 \vspace{-3mm}
\subsection{Basic graph-related definitions}\label{def_graph}
Consider a network with $N$ nodes, modeled as a static, directed graph $\mathcal{G}{=}\{\mathcal{V},\mathcal{E}\}$, where $\mathcal{V}{=}\{1,\cdots,N\}$ is the set of vertices (the nodes), and $\mathcal{E}{\subseteq}\mathcal{V}{\times}\mathcal{V}$ is the set of edges (the communication links). A directed edge from $i$ to $j$ is denoted by $(i,j){\in}\mathcal{E}$, so that information flows from $i$ to $j$. We assume $(i,i)\notin\mathcal E, \forall i \in \mathcal V$, and denote the set of {\it in-} and {\it out-neighbors} of node $i$ as
$\mathcal{N}_i^-{=}\{j:(j,i){\in}\mathcal{E}\}$ and $\mathcal{N}_i^+{=}\{j:(i,j){\in}\mathcal{E}\}$,
with cardinality $d_i^-$ ({\it in-degree}) and $d_i^+$ ({\it out-degree}), respectively. 
We will consider  strongly connected digraphs. 
\begin{definition} \label{def:graph}  A digraph $\mathcal{G}$ is strongly connected if, $\forall i,j\in \mathcal V$ with $i\neq j$, there exists a directed path from $i$ to $j$.
\end{definition}  

\begin{figure}[t]
	\centering
	\includegraphics[width = 2 in]{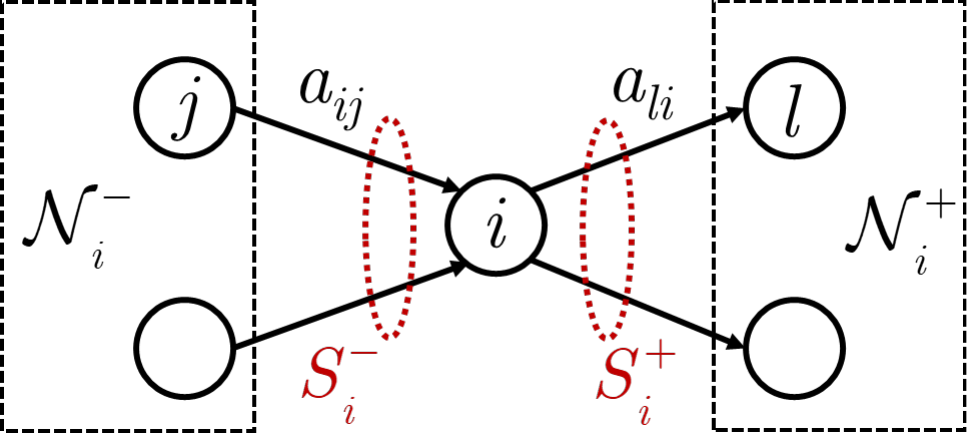}
	\caption{Some basic graph definitions.}		\label{system_model}%
	\vspace{-5mm}
\end{figure}

Associated with the digraph $\mathcal{G}$, we define a weight matrix $\mathbf{A}\triangleq (a_{ij})_{i,j=1}^N\in \mathbb{R}^{N\times N}$ such that
\begin{equation}\label{eq:Weight-matrix}\begin{cases}
a_{ij}> 0, & \text{if }(j,i)\in \mathcal{E};\\
a_{ij}= 0, & \text{otherwise};
\end{cases}   
\quad\forall i,j\in\mathcal V,
\end{equation}
along with the following quantities (cf. Fig.~\ref{system_model}) instrumental to formulate the weight-balancing problem.
\begin{definition}[In-flow, out-flow and weight-balance]  \label{def:S}
	Given a digraph $\mathcal{G}$ with weight matrix $\mathbf{A}$, 
	the {\it in}-flow of node $i$ is defined as   $S_i^-{\triangleq}\sum_{j \in \mathcal{N}_i^-}{a_{ij}}$ while the  {\it out}-flow is  $S_i^+{\triangleq}\sum_{j \in \mathcal{N}_i^+}{a_{ji}}$. The 
	  weight-balance of node $i$ is defined as $b_i{\triangleq} S_i^-{-}S_i^+$; and the overall  weight-balance vector  is ${\bf b}\triangleq(b_i)_{i=1}^N=(\mathbf A-\mathbf A^T)\mathbf 1$.
\end{definition}
\begin{definition}[Weight-balanced digraph]\label{def:WB}
A weight matrix $\mathbf A \geq {\bf 0}, \mathbf A \neq {\bf 0}$, associated to the digraph $\mathcal G$,
 is said to be \emph{weight-balanced} if
  it  induces zero balance, i.e., ${\bf b}=(\mathbf A-\mathbf A^T)\mathbf 1=\mathbf{0}$.
 \end{definition}
\vspace{-5mm}

\section{Summary of The Proposed Algorithms} \label{section_summary}
In this section, we introduce the proposed distributed quantized  weight-balancing and   consensus algorithms;  their detailed analysis will be carried out in Sec. \ref{section_weight_balance} and Sec.\ref{section_consensus}.
\vspace{-4mm}
\subsection{System model and problem formulation}
 \subsubsection*{Average consensus problem}
Let $\mathcal{V}{=}\{1,\cdots,N\}$ be a set of 
$N$ nodes. Each node $i$ controls and iteratively updates a local variable $y_i$, whose initial value is set to $y_i(0)$.  
The \emph{average consensus problem} consists in the following iterative algorithm (or variations of it): given   $\mathbf{y}(k)=(y_i(k))_{i\in \mathcal{V}}$ at time $k$, let \vspace{-0.2cm}

 \begin{equation} \label{eq:consensus_alg}
 {\bf y}(k+1) = {\bf A}{\bf y}(k),
 \end{equation}
where $\mathbf A$ is a suitably chosen weight matrix compliant with the graph [cf.\eqref{eq:Weight-matrix}]. The goal is to locally estimate the average of the initial values\vspace{-0.2cm} \begin{equation}\label{eq:consensus}
\bar{y}(0) \triangleq \frac{1}{N}\sum_{i = 1}^N{y_i(0)},
\end{equation}
i.e., $\|\mathbf{y}(k)- \bar{y}(0) \mathbf{1}\|\to 0$  as $k\to \infty$. 



We consider a setting where:    i) communications among the agents are quantized using a {\it finite} number of bits; and ii) information exchanges  flows according to the edge directions of the graph $\mathcal{G}$ (simplex communications). This puts in jeopardy the convergence of the vanilla consensus algorithm (\ref{eq:consensus_alg}), as  communications therein subsume an infinite number of bits  and    $\mathbf A$ needs to be {\it balanced} \cite{Lee2017conf}, a condition that cannot be enforced a priori without using a centralized controller with knowledge of  $\mathcal{G}$. To cope with these two issues, we first introduce   a distributed quantized weight-balancing algorithm solving the weight-balancing problem 
(cf. Sec, \ref{seb_sec_weight-balancing}); and then we integrate this algorithm with a distributed consensus algorithm using quantized simplex communications solving the average consensus problem (cf. Sec. \ref{sec:distributed_quantization}).
 \vspace{-3mm}

\subsection{Distributed quantized weight-balancing}\label{seb_sec_weight-balancing}
We propose  a  distributed, iterative algorithm to solve the weight-balancing problem over a strongly connected digraph $\mathcal{G}$ using only quantized information and simplex communications.   Note that strong connectivity guarantees the existence of a matrix, compliant to the digraph $\mathcal{G}$ (cf. D.\ref{def:graph}) that is weight-balanced (cf. D.\ref{def:WB}) \cite{Loh1970}. The proposed algorithm is formally stated in Algorithm~\ref{alg_dwb} and discussed next. 

 Each node $i$ controls the in-neighbors weights  $(a_{ij}(k))_{j \in \mathcal{N}_i^-}$.
 In \texttt{WB.1}, each node $i$ quantizes the local balance $b_i(k)$ via \eqref{signal_wb_alg_ac}, using $B_i^{\rm (w)}(k)$ bits (a $B$-bit quantizer has $2^B+1$ quantization levels), and broadcasts the quantized signal $n_i(k)$ to its out-neighbors. 
 In \texttt{WB.2}, each agent $i$ collects the signals from its in-neighbors, and updates the corresponding weights according to \eqref{a_ij}. The balance of each node is then updated according to \eqref{update_b_i}. Roughly speaking, by \eqref{a_ij}-\eqref{update_b_i} there is a {\it transfer} of the balance among nodes in the network: the quantity $\gamma(k)\,d^+_i\,n_i(k)$ (with $\gamma(k)$ denoting the step-size) is subtracted from the balance $b_i(k)$ of node $i$ [cf.~(\ref{update_b_i})], and equally divided among its out-neighbors $j{\in}\mathcal{N}_i^+$, which will increase their incoming weight $a_{ji}(k)$ by $\gamma(k)\,n_i(k)$ [cf.~\eqref{a_ij}]. Note that 1) although $n_i(k)$ may be negative, ${\bf A}(k)$ remains compliant to $\mathcal G$, which will be shown in T.\ref{thm_wb}; 2) Algorithm~\ref{alg_dwb} is fully distributed: each node $i$ only needs to know its in- and out-degrees $d_i^+$ and $d_i^-$, and to agree on a common step-size rule $\{\gamma(k)\}_{k \in \mathbb{Z}_+}$. This assumption, along with knowledge of $d_i^+$ or its equivalent information, is commonly used in distributed algorithms over directed graphs; see,  e.g., \cite{Rikos2014,Rikos2018_delay,Kempe2003,Hendrickx2015,Bernadette2018}.
 Convergence of   Algorithm~\ref{alg_dwb} is studied in Sec. \ref{section_weight_balance}.
\vspace{-0.3cm}
\subsection{Distributed quantized average consensus algorithm}\label{sec:distributed_quantization}
  We now introduce the proposed  distributed quantized average consensus algorithm over non-balanced digraphs, as described in  Algorithm \ref{alg_dac}. The algorithm  combines Algorithm~\ref{alg_dwb} with a variation of the quantized average consensus protocol based on probabilistic quantization, which we recently proposed in \cite{Lee2017conf}. 
  The  algorithm is designed so that these two building blocks run \emph{on the same time-scale}. 

  More specifically,    each node $i$ controls two set of variables, namely: i)    the in-neighbors weights  $(a_{ij}(k))_{j \in \mathcal{N}_i^-}$; and ii) the local estimate $y_i(k)$. The goal is to update these variables so that asymptotically the average consensus problem is solved while the weights converge to a balanced matrix. 
  At each iteration $k$, node $i$ quantizes its local estimate $y_i(k)$, by first clipping it within the quantization range $[q_{\min},q_{\max}]$ [cf.~(\ref{clippedY})], followed by the probabilistic quantization \eqref{quantizer} with $B_i^{\rm (c)}(k)$ bits; it then transmits the resulting quantized signal $x_i(k)$ (along with $n_i(k)$ for the weight-balancing) to its out-neighbors (\texttt{AC.1}). Upon receiving the signals $(n_j(k),x_j(k))_{j \in \mathcal{N}_i^-}$ from its in-neighbors, node $i$ updates its weights $(a_{ij}(k))_{j \in \mathcal{N}_i^-}$ using \eqref{a_ij}, and the local variable $y_i(k)$ according to \eqref{update_consensus}. The update in \eqref{update_consensus} aims at forcing a consensus on the average $\bar{y}(0)$ among the local variables $y_i(k)$. In fact, the third term in \eqref{update_consensus} is instrumental to align the local copies $y_i(k)$, while the second term $+\alpha(k)b_i(k)x_i(k)$ is a correction needed to preserve the average of the iterates, i.e., $(1/N)\sum_i y_i(k{+}1){=}(1/N)\sum_i y_i(k)$, for all $k\in \mathbb{Z}_+$ [cf.  \eqref{average_preserve}].
Hence, if all $y_i(k)$ are asymptotically  consensual, it must be $\big|y_i(k)-(1/N)\sum_i y_i(k)\big|=\big|y_i(k)-(1/N)\sum_i y_i(0)\big| \underset{k\to\infty}{\longrightarrow} 0$. 

Convergence of Algorithm \ref{alg_dac} is studied in Sec. \ref{section_consensus}.

\begin{algorithm}[t]
\caption{Distributed Quantized Weight-Balancing} 
\label{alg_dwb} 
\begin{algorithmic} 
	\Require  \texttt{(WB.0)} ${\bf A}(0)$; $\big\{\gamma(k),\big(B_i^{\rm (w)}(k)\big)_{i \in \mathcal{V}}\big\}_{k \in \mathbb{Z}_+}$. \\
	Set $k=0$. Repeat \texttt{(WB.1)-(WB.2)} for $k=1,2,\dots$ until a termination criterion is satisfied;
    \State \texttt{(WB.1)} Each node $i$ broadcasts the $n_i(k)$ to $\mathcal{N}_i^+$, where 
    \begin{align}
    \!\!\!n_i(k) &{=} 
    {\rm sgn}(b_i(k))
    \!\min\!\left\{\!\!\frac{2}{2^{B_i^{\rm (w)}(k)}}\!\!\left[\left\lceil\!\frac{2^{B_i^{\rm (w)}(k)}\vert b_i(k) \vert}{2d_i^+\gamma(k)}\!\right\rceil \!\!{-} 1\right]\!\!, 1\right\}. \label{signal_wb_alg_ac}
\end{align}
    \State \texttt{(WB.2)} Each node $i$ collects $n_j(k)$ from its in-neighbors $j \in \mathcal{N}_i^-$, and updates 
    \begin{align}
    a_{ij}(k+1) &= a_{ij}(k)+\gamma(k)\, n_j(k), \qquad\forall j \in \mathcal{N}_i^-, \label{a_ij}\bigskip \\
    \!\!b_i(k+1) &= b_i(k)-\gamma(k)\,d_i^+\,n_i(k)\!+\!\gamma(k)\!\!\sum_{j\in \mathcal{N}_i^-}\!\!n_j(k).\label{update_b_i}
    \end{align}
\end{algorithmic}
\end{algorithm}

\begin{algorithm}[t]
\caption{Distributed Quantized Average Consensus} 
\label{alg_dac} 
\begin{algorithmic} 
	\Require  \texttt{(AC.0)} Init. Algorithm \ref{alg_dwb}  as in \texttt{(WB.0)}; $q_{\min}, q_{\max}$; $\big\{\alpha(k),\big(B_i^{\rm (c)}(k)\big)_{i \in \mathcal{V}}\big\}_{k \in \mathbb{Z}_+}$; and ${\bf y}(0)$.\\
	Set $k=0$. Repeat \texttt{(AC.1)-(AC.3)} for $k=1,2,\dots$ until a termination criterion is satisfied;
    \State \texttt{(AC.1)} Each node $i$ broadcasts $n_i(k)$ (cf. \eqref{signal_wb_alg_ac}) and $x_i(k)$ to $\mathcal{N}_i^+$, where $x_{i}(k)=0$ if $B_i^{\rm (c)}(k) =0$ and, if $B_i^{\rm (c)}(k)>0$,
\begin{align} 
    &\!\!\!\!\!\!\!\!x_{i}(k) \!=\!\left\{	\begin{array}{*{20}l} \!\!\!\!q_{\min}{+}\!\left\lceil\! \dfrac{\tilde{y}_i(k)\!-\!q_{\min}}{\Delta\left(B_i^{\rm (c)}(k)\right)}\!\right\rceil \!\Delta\!\left(\!B_i^{\rm (c)}(k)\!\right), \text{ w.p. } p_i(k); \medskip \\
	\!\!\!\!q_{\min}{+}\!\left\lfloor\! \dfrac{\tilde{y}_i(k)\!-\!q_{\min}}{\Delta\left(B_i^{\rm (c)}(k)\right)}\!\right\rfloor \!\Delta\!\left(\!B_i^{\rm (c)}(k)\!\right), \text{ w.p. } 1\!-\!p_i(k),
	\end{array}\right.\!\!\!\!\!\!\!\!\!\!\!\!\!\label{quantizer}
    \end{align}
    where $\tilde{y}_i(k) = {\rm clip}(y_{i}(k); q_{\min}, q_{\max})$,
    \begin{align}
\Delta(B) &= \frac{q_{\max}-q_{\min}}{2^B-1}, \forall B \in \mathbb{Z}_{++}, \label{clippedY} \\
p_i(k) &= \frac{\tilde{y}_i(k)-q_{\min}}{\Delta\left(B_i^{\rm (c)}(k)\right)} - \left\lfloor \frac{\tilde{y}_i(k)-q_{\min}}{\Delta\left(B_i^{\rm (c)}(k)\right)}\right\rfloor. \nonumber
	\end{align}
    \State \texttt{(AC.2)} Each node $i$ collects $(n_j(k),x_j(k))$ from its in-neighbors $j{\in}\mathcal{N}_i^-$, updates $\left(a_{ij}(k+1)\right)_{j{\in}\mathcal{N}_i^-}$ (cf. \eqref{a_ij}) and
    \begin{align}
    \!\!y_i(k+1) & = y_i(k)+\alpha(k)b_i(k)\,x_i(k) \nonumber\\
     &\quad +\alpha(k)\!\!\!\sum\limits_{j \in \mathcal{N}_i^-}{\!\!\!a_{ij}(k)\,\Big(x_j(k)-x_i(k)\Big)}.  \label{update_consensus}
    \end{align}
\end{algorithmic}
\end{algorithm}
\vspace{-1mm}

\section{Distributed Quantized Weight-Balancing} \label{section_weight_balance}
We study  convergence   of Algorithm \ref{alg_dwb} under the following mild assumptions.\footnote{The analysis can be extended to the case in which each node uses its own step-size $\{\gamma_i(k)\}$, provided that: 1) every node knows the step-size of its in-neighbors, and 2) every $\{\gamma_i(k)\}$ satisfies A.\ref{assump_gamma2}.} \vspace{-2mm}

\begin{assumption} \label{assump_Bw}
Let $\left\{B^{\rm (w)}(k)\right\}_{k \in \mathbb{Z}_+}$ be a sequence satisfying $B^{\rm (w)}(k){\in}\{0,1\}$ and $ \sum_{t = nW}^{(n+1)W-1}{B^{\rm (w)}(t)} \geq 1,$ for all $k,n \in \mathbb{Z}_+$ and some  $W\in \mathbb{Z}_{++}$. Then, there exists $B_{\max} \in \mathbb Z_+$ such that the number of bits $\big\{B_i^{\rm (w)}(k)\big\}_{k \in \mathbb{Z}_+}$ satisfies: for all  $i\in \mathcal{V}$,
\begin{align*}
\begin{cases}
B_{\max}\geq B_i^{\rm (w)}(k)\geq B^{\rm (w)}(k), & \text{if }B^{\rm (w)}(k)=1; \\
B_i^{\rm (w)}(k)= 0, & \text{else.}
\end{cases}
\end{align*}  \vspace{-4mm}
\end{assumption}

\begin{assumption} \label{assump_gamma2}
The step-size $\{\gamma(k)\}_{k \in \mathbb{Z}_+}$ and initial weight matrix $\mathbf{A}(0)\triangleq (a_{ij}(0))_{i,j=1}^N$ satisfy:
\begin{align}
&\gamma(k)  = c_1^{-n},\text{ with } n{\in}\mathbb{Z}_+ \!:\! (c_1^n \!-\! 1)c_2  \leq  k  \leq  (c_1^{n+1} \!-\! 1)c_2 \!-\! 1; \nonumber \\
&\text{and }\quad a_{ij}(0) = \mathcal I\{(j,i)\in \mathcal E\}, \label{a_ij0}
\end{align}
respectively, where $c_1\in\mathbb Z,c_1\geq 2$, and $c_2 \in \mathbb R_{++}$.

 Define $\bar{\gamma}(k) \triangleq \gamma(k)2^{1-B_{\max}^{(w)}}$, $k\in \mathbb{Z}_+$.
\end{assumption}\vspace{-0.1cm}
 Note that the step-size satisfying A.\ref{assump_gamma2}  is vanishing and non-summable, as shown below.\vspace{-0.2cm}
 \begin{lemma} \label{lemma_gamma_ub}
If $\{\gamma(k)\}_{k \in \mathbb{Z}_+}$ satisfies A.\ref{assump_gamma2}, then
\begin{align}
\frac{c_2}{k+c_2} \leq \gamma(k) \leq \frac{c_1\,c_2}{k+c_2}, \quad \forall k\in \mathbb{Z}_+.\label{eq:lower_upper_gamma}
\end{align}
\end{lemma}
\begin{IEEEproof}
Let $n\!\in\! \mathbb{Z}_+$.
Note that 
$\gamma(k){=}c_1^{-n},\forall k:(c_1^n\!-\!1)c_2 \!\leq\! k \!\leq\! (c_1^{n+1}\!-\!1)c_2\!-\!1$. Then, the upper and lower bounds on
$\gamma(k)$ are obtained by bounding $c_1^n$ with respect to this interval.\hfill
\end{IEEEproof}
A.\ref{assump_gamma2}  is consistent with similar choices adopted  in stochastic optimization \cite{Bertsekas1989}, such as   $\gamma(k){=}1/(k{+}1)$. However the diminishing and non-summability  properties alone are not sufficient to prove convergence of Algorithm~\ref{alg_dwb};     A.\ref{assump_gamma2} further guarantees that the balance at each node is always an integer multiple of the current step-size (L.\ref{lemma_eps_gamma}, cf. App. \ref{append:lemma_decr_event}
), which will be shown to be a key property to prove that $\Vert {\bf b}(k)\Vert_1$ is asymptotically vanishing. 
An instance  of $\{\gamma(k)\}_{k \in \mathbb{Z}_+}$ satisfying A.\ref{assump_gamma2} is \cite{Lee_2018CDC}\vspace{-0.1cm}
\begin{align}
\{\gamma(k)\}_{k \in \mathbb{Z}_+} = \left\{1,\frac{1}{2},\frac{1}{2},\frac{1}{4},\frac{1}{4},\frac{1}{4},\frac{1}{4},\frac{1}{8},\cdots \right\}. \label{gamma_ex}
\end{align}
Note that a fixed step-size $\gamma(k) = \gamma$, for all $k$, may fail to achieve convergence. In fact, it is possible that each $0<\vert b_i(k) \vert \leq d_i^+\gamma$, so that $n_i(k) = 0$ $\forall i$   and there is no further transfer of balance, resulting still in an unbalanced digraph. 

A.\ref{assump_Bw} states that at least once over a time window of duration $W$, \emph{all} nodes are simultaneously communicating at least one bit to their out-neighbors. This offers some flexibility in the design of the communication protocol. For instance, nodes can transmit one bit at each time slot \cite{Lee_2018CDC}, yielding one bit per channel use, or transmit one bit every $W>1$ time slots, resulting in a lower effective rate of $1/W$ bits per channel use. \vspace{-0.4cm}
\\\indent We are now ready to state our main convergence result. \vspace{-0.2cm} \begin{theorem}\label{thm_wb}
Let $\{\mathbf{A}(k)\}_{k\in \mathbb{Z}_+}$ be the sequence generated by Algorithm~\ref{alg_dwb}
  under A.\ref{assump_Bw} and A.\ref{assump_gamma2}.
   Then, there hold: 
   \begin{align*}\vspace{-0.2cm}
&\text{(a) }\Vert {\bf b}(k)\Vert_1 = \mathcal{O}\left(\frac{1}{k}\right),
\ \ \ \  \text{(b) }\lim\limits_{k\rightarrow \infty}{{\bf A}(k)} = {\bf A}^*,
\\
&\text{(c) }0 < a_{\min} \leq a_{ij}(k) \leq a_{\max}, \forall (j,i) \in \mathcal E, \forall k \in\mathbb Z_+,
   \end{align*}
   where ${\bf A}^*$ is weight-balanced and $a_{\min}$, $a_{\max}$ are defined in \eqref{a_globalbound}.
\end{theorem}\vspace{-0.5cm} 

\subsection{Proof of Theorem~\ref{thm_wb}} \label{subsec_wb_proof}
\subsubsection{Proof of statement (a)}
\begin{figure}[t!]\vspace{-0.2cm}
	\centering
	\includegraphics[width = \linewidth]{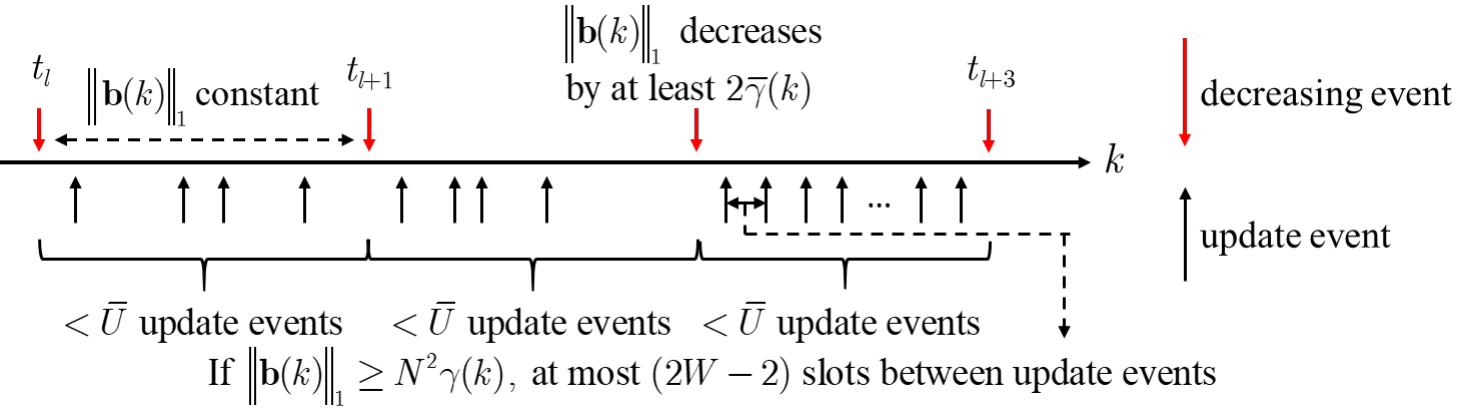} 
	\caption{Key properties of $\Vert {\bf b}(k)\Vert_1$ in Algorithm \ref{alg_dwb}.
	}\label{thm1_properties}\vspace{-0.5cm}
\end{figure}

     We begin by highlighting the main steps of the proof,   with the help of Fig. \ref{thm1_properties}.  {\bf{Step 1}:} 
     We show that $\{\Vert {\bf b}(k)\Vert_1\}_{k\in \mathbb{Z}_+}$ is non-increasing. Furthermore, we identify two key events affecting the dynamics of $\{\Vert {\bf b}(k)\Vert_1\}_{k\in \mathbb{Z}_+}$,
      namely: the  so-called ``decreasing event'' $\mathcal{D}_k$ and ``update event'' $\mathcal{U}_k$.     $\mathcal{D}_k$, formally defined in \eqref{event_D}, occurs if at timeslot $k$ either one of the following two facts happen: 1) a node 
      transfers its nonzero balance to an out-neighbor with balance of opposite sign;   2) two nodes, with balance of opposite sign,  
      transfer their balances to a common out-neighbor. 
     On the other hand,  $\mathcal{U}_k$, formally defined in  \eqref{event_U},  occurs if at timeslot $k$ a node transfers its balance to its out-neighbors (i.e., $n_i(k) {\neq}0$,  for some $i\in\mathcal{V}$) but $\mathcal{D}_k$ does {\it not} occur.   Note that 
     it can happen that  neither $\mathcal{D}_k$ nor $\mathcal{U}_k$ occur at some $k$; this is the case when $\vert b_i(k) \vert$ is ``too small'' or all nodes are inactive ($B^{\rm (w)}(k)=0$). 
     We show  that $\Vert {\bf b}(k)\Vert_1$ decreases by at least $2\bar\gamma(k)$, with  $\bar{\gamma}(k) \triangleq \gamma(k)2^{1-B_{\max}^{(w)}}$, iff. $\mathcal D_k$ occurs, and remains unchanged otherwise (cf. L.\ref{lemma_decr_event}). {\bf{Step 2:} }
      To guarantee that $\{\Vert {\bf b}(k)\Vert_1\}_{k\in \mathbb{Z}_+}$ vanishes, the decreasing event  must occur \emph{sufficiently often}. Towards this end, we prove two key properties of the decreasing and update events, namely: 
      
\begin{description}
	\item[{\bf P1)}]   there are at most $\bar{U} \triangleq N^{2N-2}$ update events between two consecutive decreasing events;
   \item[{\bf P2)}]  
    if $\Vert {\bf b}(k)\Vert_1{\geq}N^2\gamma(k)$ (roughly speaking, if $\{\Vert {\bf b}(k)\Vert_1\}_{k\in \mathbb{Z}_+}$ does not decrease sufficiently fast),   there are at most $2W-2$ timeslots between two consecutive update events.\end{description}    
    The decreasing step-size together with P2 guarantee that update events occur within bounded time; this property combined with P1 guarantees that decreasing events 
   occur at uniformly bounded time intervals. Finally, \textbf{Step 3} builds on the above results to prove statement (a) of the theorem. Roughly speaking,   
   one can infer that:
     either 1) $\{\Vert {\bf b}(k)\Vert_1\}_{k\in \mathbb{Z}_+}$ is below the diminishing  threshold (Step 1), causing it to vanish (since the step-size vanishes); or 2) it exceeds the threshold at some timeslots, 
    causing it to be suppressed by the decreasing event 
    (Step 2) until it falls again below the vanishing threshold. 
    
    We proceed next with the  formal proof.

\noindent\textbf{Step 1:} We begin by introducing the   definition of    $\mathcal D_k$ and 
$\mathcal U_k$.
\begin{definition}[Decreasing event $\mathcal D_k$ and its occurrence time $t_l$] 
Let  $\mathcal D_k$,   $k\in\mathbb{Z}_+$,    be defined  as 
	 	\begin{align}
	 \!\!\!\!\!\exists i, (j_1,j_2)\in [\mathcal N_i^-]^2: b_i(k)n_{j_1}(k){<} 0\vee n_{j_2}(k)n_{j_1}(k){<} 0
	 .	\label{event_D}
	 \end{align}
	  \label{def_event_D} 
Furthermore, let $t_l$, $l\in\mathbb Z_+$, be the timeslot of occurrence of the $l$th decreasing event; recursively, $t_{0} \triangleq -1$ and  $l \in \mathbb Z_{++}$,
\begin{align}
t_l\triangleq\min\{k>t_{l-1}:\mathcal I\{\mathcal D_k\}=1\},\ (\text{possibly, }t_l=\infty). \label{t_l}
\end{align}
\end{definition}
\begin{definition}[Update event $\mathcal U_k$]
Let $\mathcal U_k$,  $k\in\mathbb{Z}_+$, be defined as:\vspace{-0.2cm}
\begin{equation}
\label{event_U}
\exists\   i  \in \mathcal V\,:\, n_i(k) \neq 0 \quad \wedge\quad  \mathcal{I}(\mathcal D_k)= 0.
\end{equation}
\end{definition}

We show next that $\{\Vert {\bf b}(k)\Vert_1\}_{k\in \mathbb{Z}_+}$ decreases by at least $2\bar{\gamma}(k)$ iff. $\mathcal D_k$ occurs, and remains unchanged otherwise.
\begin{lemma} \label{lemma_decr_event} There holds
 \begin{align}\label{eq:lemma3}
 	\begin{cases}
\Vert{\bf b}(k+1)\Vert_{1}\leq \Vert{\bf b}(k)\Vert_{1}-2\,\bar{\gamma}(k), & \text{if }\mathcal{I}\left(\mathcal D_k\right)=1;\\
\Vert{\bf b}(k+1)\Vert_{1}= \Vert{\bf b}(k)\Vert_{1}, & \text{otherwise.}
\end{cases}
 \end{align}
\end{lemma}
\begin{IEEEproof}
See App.A-\rom{2}.\smallskip 
\end{IEEEproof}

\noindent\textbf{Step 2:}
This step characterizes ``how often'' $\mathcal D_k$ and $\mathcal U_k$ occur (properties P1 and P2), and the implication on $\Vert {\bf b}(k)\Vert_1$. We first provide some intuition motivating our approach.
 
 \noindent   {\bf $\bullet$ Intuition:} Let us look at the balance transfer within two consecutive decreasing events at (finite) times $t_l$ and $t_{l+1}$. Let 
 \begin{align}
 \label{V+V-}
 \mathcal V^+(k) \triangleq \{i:b_i(k) > 0\},\ \mathcal V^-(k) \triangleq \{i:b_i(k) < 0\}
 \end{align}
be the set of nodes with positive and negative balance. Note that  $\mathcal V^+(k), \mathcal V^-(k)\neq \emptyset$  iff. $\Vert {\bf b}(k)\Vert_1 > 0$, since $\sum_{i=1}^N b_i(k)=0$. 
For the next decreasing event $\mathcal{D}_{t_{l+1}}$ to occur at time $t_{l+1}$: either 1) a node $i{\in}\mathcal{V}^+(t_{l+1})$ (resp. $i{\in}\mathcal{V}^-(t_{l+1})$) has enough balance (i.e., $n_i(t_{l+1})\neq 0$) to trigger the update to an out-neighbor in $\mathcal{V}^-(t_{l+1})$ (resp. $\mathcal{V}^+(t_{l+1})$); 
or 2) two nodes $j_1\in\mathcal{V}^+(t_{l+1})$ and $j_2\in\mathcal{V}^-(t_{l+1})$   trigger an update to an out-neighbor in     $\mathcal{N}_{j_1}^+\bigcap \mathcal{N}_{j_2}^+$ $(\neq \emptyset)$. This balance is built-up throughout  the update events within $(t_l,t_{l+1})$, during which nodes in $\mathcal{V}^+(k)$ (resp. $\mathcal{V}^-(k)$), $k\in (t_l,t_{l+1})$, keep transferring  part of their balance towards the  out-neighbors in $\mathcal{V}^+(k)$ (resp. $\mathcal{V}^-(k)$) that are {\it closer} to   nodes outside $\mathcal{V}^+(k)$  (resp. $\mathcal{V}^-(k)$). Hence, one can expect that the decreasing event $\mathcal{D}_{t_{l+1}}$ will occur after a certain number of update events, specifically when a sufficient amount of balance has been transferred to some nodes having out-neighbors in $\mathcal V \setminus \mathcal{V}^+(k)$ (resp. $\mathcal V \setminus \mathcal{V}^-(k)$). To characterize this number and show that it is bounded over each interval $(t_l,t_{l+1})$, the proposed idea is to construct a {\it nonnegative, integer-valued} function of $\mathbf b(k),k{\in}(t_l,t_{l+1})$, denoted by $U(k)$, which (a) {\it strictly increases} whenever $\mathcal U_k$ occurs; and (b) is  uniformly upper {\it bounded} on $(t_l,t_{l+1})$. These properties   guarantee that the number of update events within $(t_l,t_{l+1})$ is bounded, which proves P1. The same function $U$ will   be also used to prove   P2 (cf. P.\ref{lemma_Uk} \& C.\ref{corollary_minimal_decrease}). Next, we build $U(k)$ and prove P1 and P2.  
 
\noindent {\bf $\bullet$ Building the function $U(k)$:} Let $l \in \mathbb Z_+$, $t_l < \infty$,
\begin{equation}
U_h \triangleq N^{2(N-1-h)},  \quad h\in\{1,2,\dots,N-1\}. \label{U_n}
\end{equation}
We define the set (possibly empty) of nodes that are $h$ directed hops away from a node with opposite sign of balance as
\begin{align}
\mathcal{V}_h(k) \triangleq&  \left\{ i\in \mathcal{V}^+(k)\,:\,\underset{j \in \mathcal{V}^{-}(k)}\min d(i,j) = h\right\} \nonumber \\
&\!\!\!\!\!\!\cup \left\{ i\in \mathcal{V}^-(k)\,:\,\underset{j \in \mathcal{V}^{+}(k)}\min d(i,j) = h\right\},
\end{align}
where $d(i,j)$ is the directed distance between $i$ and $j{\in}\mathcal V$.
Then, we define the function
\begin{align}
U(k) {\triangleq}&  \sum\limits_{h = 1}^{N-1}{U_h\!\!\!\!\sum\limits_{i \in \mathcal{V}_h(k)}{\!\!\!\min\left\{\frac{\vert b_i(k) \vert}{\bar\gamma(k)}, N\right\}}},\,\,\,  k{\in}(t_l, t_{l+1}].\label{U_k} 
\end{align}
Based on the above discussion, nodes in $\mathcal V_1(k)$   are   ``more important''   than  nodes in $\mathcal{V}_2(k)$, in the sense that they will more immediately trigger the next decreasing event;   nodes in  $\mathcal{V}_2(k)$ are more important than those in $\mathcal{V}_3(k)$, and so on. The   function $U(k)$ aims at capturing this hierarchical   transfer of balance  along the chain  $\mathcal{V}_{N-1}(k) \rightarrow \cdots \mathcal{V}_{h}(k)\cdots \rightarrow \mathcal{V}_{1}(k)$ during each update event, up until $\mathcal D_{t_{l+1}}$ occurs. In particular, we want $U(k)$ to increase     its value by (at least) one (integer) unit every time one of such transfers happens  (i.e., $\mathcal U_{k}$ occurs).  

To motivate the choice of \eqref{U_k}, let us look at the balance transfer during an update event $\mathcal U_{k}$ at time $k{\in}(t_l,t_{l+1})$; for the sake of simplicity, say $\mathcal U_{k}$ is triggered by node $i{\in}\mathcal{V}_{h+1}(k)$. As a result, node $i$ transfers  part of its balance to its out-neighbors $\mathcal N_i^{+}{\cap}\mathcal{V}_h(k)$,\footnote{Note that $\mathcal N_i^{+}\cap\mathcal{V}_h(k)\neq\emptyset$ as, by definition of $\mathcal V_h(k)$, at least one node in $\mathcal N_i^{+}$ is one step closer to nodes with balance of opposite sign.} according to \eqref{update_b_i}. In \eqref{update_b_i}, $\bar\gamma(k)$ can be regarded as the unit of balance and  $\vert b_i(k)\vert/\bar\gamma(k)$ is the normalized imbalance, an integer number (cf. L.\ref{lemma_eps_gamma}, App.A-\rom{3}).  Such a node $i{\in}\mathcal{V}_{h+1}(k)$ experiences a decrease of its normalized imbalance by $d_i^+$ units while the normalized imbalance of $j{\in}\mathcal{V}_{h}(k)$ increases by at least one unit. To encode this balance transfer as an increase of $U(k)$ by at least one integer, we can associate  it with the ``carry on'' operation from a digit to the next more significant one in a positional notational  representation of the $U(k)$ value.   Specifically, $U(k)$ is   expressed in radix-$N^2$ notation wherein the sum normalized imbalance of nodes in $\mathcal{V}_1$ contributes to the most significant digit, the one of nodes in $\mathcal{V}_2$ contributes to the second most significant digit, and so on. By doing so, when $\mathcal U_k$ occurs as above, the aforementioned  exchange of  balance  $\mathcal{V}_{h+1}(k)\rightarrow \mathcal{V}_h(k)$ triggers the transfer of one unit from the   $(h+1)${\it th} most representative digit to the $h${\it th} one, so that $U(k)$ increases by at least one unit.

\noindent {\bf $\bullet$ Proof of P1 and P2:}  P.\ref{lemma_Uk} below  states the desired properties of $U(k)$ and proves P2 as a by-product; P1 follows from C.\ref{corollary_minimal_decrease}. 
\begin{proposition}[Properties of $U(k)$] \label{lemma_Uk}
Let $k{\in}(t_l, t_{l+1})$. Then:
\begin{enumerate}[(i)]
\item $U(k)$ is non-negative, non-decreasing, and 
upper bounded by $U(k) \leq N^{2N-1}$;
\item $U(k)$ strictly increases if $\mathcal U_k$ occurs, i.e.,   $
  U(k+1)\geq U(k)+N\mathcal{I}(\mathcal{U}_{k});
$
\item If $\Vert {\bf b}(k)\Vert_1\geq N^2\gamma(k)$, 
then an update or decreasing event occurs within the next $2W-1$ slots; hence,\\ $U(k+2W-1) \geq U(k)+N,\forall k\leq t_{l+1}-2W{+}1$.\footnote{Since we are interested in the variations of $U(k)$ within   $(t_l, t_{l+1}]$, the case $k>t_{l+1}-2W{+}1$ is irrelevant.}
\end{enumerate}
\end{proposition}
\begin{IEEEproof}
See Appendix~A-\rom{3}.
\end{IEEEproof}

\begin{corollary} \label{corollary_minimal_decrease}  
\begin{enumerate}[(i)]
\item There are at most $\bar{U}=N^{2N-2}$ update events between two consecutive decreasing events.
\item If $\Vert{\bf b}(k)\Vert_1 \geq N^2\gamma(k)$, $k\in\mathbb{Z}_+$, then
\begin{align*}
\left\Vert {\bf b}\!\left(k\!+\!(2W\!-\!1)\bar{U}\right)\right\Vert_1
\!\leq\!\left\Vert {\bf b}\!\left(k\right)\right\Vert_1{-}2\bar{\gamma}\left(k{+}(2W{-}1)\bar{U}\right).
\end{align*}
\end{enumerate}
\end{corollary}
\begin{IEEEproof}
{\it (i)} is a direct result of P.\ref{lemma_Uk}{(i)-(ii)} and the fact that $U(k) \geq 2$ iff. ${\bf b}(k) \neq {\bf 0}$.   {\it(ii)}: 
let $k$ such that $\Vert {\bf b}(k)\Vert_1 {\geq} N^2\gamma(k)$, and let $t_{l+1}=\{\tau{\geq}k:\mathcal I\{\mathcal D_\tau\}{=}1\}$ be the next decreasing event at or after $k$  (possibly, $t_l{=}-1$ and/or $t_{l+1}=\infty$). Invoking P.\ref{lemma_Uk}, we infer that i) there are at most $\bar{U}$ update events in $(t_l,t_{l+1})$;  ii) $\Vert {\bf b}(\tau)\Vert_1{\geq}N^2\gamma(\tau),\forall\tau {\in}(t_{l},t_{l+1}]$ ($\gamma(\tau)$ is non-increasing and $\Vert {\bf b}(\tau)\Vert_1$ only decreases after decreasing events), so that the first update or decreasing event after $t_l$ occurs within $2W-1$ timeslots and subsequent ones are separated by at most $2W-2$ timeslots until the next decreasing event at $t_{l+1}$. These two facts together imply $t_{l+1}\leq k+(2W-1)\bar{U}-1$;  therefore, 
\begin{align*}
&\left\Vert {\bf b}\!\left(k\!+\!(2W\!-\!1)\bar{U}\right)\right\Vert_1
\!\overset{\eqref{eq:lemma3}}{\leq}\! \left\Vert {\bf b}\!\left(t_{l+1}\!+\!1\right)\right\Vert_1 \\
&\!\overset{\eqref{eq:lemma3}}{\leq}\!\left\Vert {\bf b}\!\left(t_{l+1}\right)\right\Vert_1\!-\!2\bar{\gamma}\!\left(t_{l+1}\right)
\!\!\!\overset{(A.\ref{assump_gamma2})}{\leq}\!\!\left\Vert {\bf b}\!\left(k\right)\right\Vert_1{-}2\bar{\gamma}\!\left(k\!+\!(2W\!-\!1)\bar{U}\right).
\end{align*}
\end{IEEEproof}\smallskip
\noindent\textbf{Step 3}: 
We now prove that $\Vert {\bf b}(k)\Vert_1{=}\mathcal O(1/k)$. Equivalently,
\begin{align} 
\exists\ 0\!<\!M\!<\!\infty \,\text{ and }\bar{k}\in \mathbb{Z}_{++} \!:\!   k\cdot\Vert {\bf b}(k)\Vert_1\!\leq\!M, \forall k\geq \bar{k}. \label{limit2}
\end{align}
To this end, note that it suffices to show that
\begin{align}
\label{nhgrug}
\Vert {\bf b}(k)\Vert_1\leq N^2c_1\gamma(k),\quad\forall k\geq \bar k.
\end{align}
In fact, using L.\ref{lemma_gamma_ub} in App.\ref{append:wb}-\rom{1},
\eqref{nhgrug} implies
$$
k\Vert {\bf b}(k)\Vert_1
\leq \frac{N^2c_1^2c_2k}{k+c_2}=\mathcal{O}(N^2c_1^2c_2),
$$
so that \eqref{limit2} readily follows with
$M\triangleq N^2c_1^2c_2$.
To prove \eqref{nhgrug}, 
let $\tilde{k}_m\triangleq(c_1^m-1)c_2,m\in\mathbb Z_+$;
we define $\bar k$ as
 \begin{align}
\bar k=\min\left\{\tau\geq \tilde k_n:\Vert {\bf b}(\tau)\Vert_1 < N^2\gamma(\tau)\right\},\label{k_bar2}
 \end{align} 
 for some sufficiently large $n{\in}\mathbb{Z}_+$ to be determined.
The existence of such $\bar{k}$  is guaranteed by L.\ref{lemma_T-} in App.\ref{append:wb}-\rom{1}. 
Let
$$p=\min\{m>n:\tilde{k}_m>\bar k\}.$$
Then, it readily follows that, for all $k\in[\bar k,\tilde{k}_p-1]$,
\begin{align}
\nonumber
  \Vert {\bf b}(k)\Vert_1 & {\overset{\eqref{eq:lemma3}}\leq}\left\Vert {\bf b}\left(\bar k\right)\right\Vert_1
{\overset{\eqref{k_bar2}}<}N^2\gamma\left(\bar k\right)  
{\overset{(A.\ref{assump_gamma2})}=}
N^2\gamma(k),
\end{align}
so that \eqref{nhgrug} holds for $k{\in}[\bar k,\tilde{k}_p-1]$.
It remains to prove that it holds for $k{\geq}\tilde{k}_p$.
We do so by induction. Assume that
it holds at $k{\in}[\bar k,\tilde{k}_m-1]$, for some $m{\geq}p$, and that 
\begin{align}
\label{gds}
\Vert {\bf b}(\tilde{k}_m-1)\Vert_1\leq N^2\gamma(\tilde{k}_m-1).
\end{align}
Clearly, this is true for $m=p$. We show next that this condition implies that
\eqref{nhgrug} holds for $k\in[\bar k,\tilde{k}_{m+1}-1]$, and  
\begin{align}
\label{gds2}
\Vert {\bf b}(\tilde{k}_{m+1}-1)\Vert_1\leq N^2\gamma(\tilde{k}_{m+1}-1).
\end{align}
Therefore, \eqref{nhgrug} holds $\forall k\geq \bar k$.
To show the induction step, note that \eqref{gds} implies
 \eqref{nhgrug},
  $\forall k\in[\tilde{k}_m,\tilde{k}_{m+1}-1]$, since
\begin{align}
\Vert {\bf b}(k)\Vert_1
&\stackrel{(L.2)}{\leq} \Vert {\bf b}(\tilde{k}_m-1)\Vert_1
\stackrel{\eqref{gds}}{\leq} N^2\gamma(\tilde{k}_m-1) \nonumber\\
&\stackrel{(A.1)}{=}N^2c_1\gamma(k)
,\ \forall k\in[\tilde{k}_m,\tilde{k}_{m+1}-1].
\label{fgh}
\end{align}
It  remains to prove \eqref{gds2}; we do it  by contradiction.
Assume \eqref{gds} and \eqref{fgh} hold but \eqref{gds2} does not. 
Then,
\begin{align*}
&N^2c_1\gamma(k)
\stackrel{\eqref{fgh}}{\geq} \Vert {\bf b}(k)\Vert_1
\stackrel{(L.2)}{\geq} \Vert {\bf b}(\tilde{k}_{m+1}-1)\Vert_1 \\
&>N^2\gamma(\tilde{k}_{m+1}-1)
\stackrel{(A.1)}{=}N^2\gamma(k)
,\ \forall k\in[\tilde{k}_m,\tilde{k}_{m+1}-1].
\end{align*}
Choosing $n$ large enough so that, for $T\triangleq 2^{B_{\max}{-}2}N^2(c_1{-}1)$, \vspace{-0.2cm}
\begin{equation}
\tilde{k}_m+T(2W-1)\bar{U}
\leq\tilde{k}_{m+1}-1,\ \forall m\geq p>n, \label{cond_T}
\end{equation}
(this is possible since
$\tilde{k}_{m+1}-\tilde{k}_m
\geq \tilde{k}_{n+1}-\tilde{k}_n
=c_2(c_1-1)c_1^n$)
we can then apply C.\ref{corollary_minimal_decrease}{\it(ii)} recursively
$T$ times, yielding
\begin{align*}
&\Vert {\bf b}(\tilde{k}_{m+1}-1)\Vert_1
\stackrel{\eqref{cond_T},(L.2)}{\leq}
\Vert {\bf b}(\tilde{k}_m+T(2W-1)\bar{U}\Vert_1 \\
&{\leq}\Vert {\bf b}(\tilde{k}_m)\Vert_1
{-}\sum_{j=0}^{T-1}2\bar{\gamma}(\tilde{k}_m{+}j(2W{-}1)\bar{U})
\stackrel{(a)}{\leq}
N^2 \gamma(\tilde{k}_{m+1}-1),
\end{align*}
\noindent where $(a)$ follows from A.\ref{assump_gamma2} and \eqref{fgh}. This proves the contradiction, hence $\Vert {\bf b}(k)\Vert_1{=}\mathcal O(1/k)$.

\subsubsection{Proof of statement (b)} 
Convergence of $\{{\bf A}(k)\}_{k \in \mathbb{Z}_+}$ to $\bf A^*$   is a consequence of the following Lemma. 
\begin{lemma} \label{lemma_a_bound}
The sequence $\{{\bf A}(k)\}_{k \in \mathbb{Z}_+}$ is a Cauchy sequence.
\end{lemma}
\begin{IEEEproof} See Appendix \ref{sec_proof_lemma_A}.\end{IEEEproof}


\subsubsection{Proof of statement (c)} 
First, using the fact that for any $x \in \mathbb R_+$, there exists $\beta \in [0, 1)$ such that $\lceil x\rceil -1 = \beta x$, it follows that
\begin{align}
n_i(k) = \beta_i(k) \frac{b_i(k)}{d_i^+\gamma(k)}, \text{ for some }\beta_i(k) {\in} [0, 1). \label{lb_1}
\end{align}
Note that \eqref{a_ij0} and \eqref{a_ij} imply $a_{ji}(k), \forall j{\in} \mathcal N_i^+$ have the same value. Let ${\bf S}^+(0) {=} {\rm diag}\{S^+(k), i {\in} \mathcal V\}, \pmb\beta(k) = {\rm diag}\{\beta_i(k), i \in \mathcal V\}, {\bf a}(k) = (a_i(k))_{i=1}^N$, where $a_i(k) = a_{ji}(k), j{\in} \mathcal N_i^+$. Applying \eqref{lb_1} to the update of $a_{ij}(k)$, it follows that
\begin{align*}
{\bf a}(k+1) = {\bf a}(k)+\pmb\beta(k)\left({\bf S}^+(0)\right)^{-1}\left[{\bf A}(0)-{\bf S}^+(0)\right]{\bf a}(k),
\end{align*}
where we use the fact ${\bf b}(k) = {\bf S}^-(k) - {\bf S}^+(k) = {\bf A}(0){\bf a}(k) - {\bf S}^+(0){\bf a}(k)$,
which implies that
\begin{align*}
{\bf a}(k+1) = \left[{\bf I}-\pmb\beta(k)+\pmb\beta(k)\left({\bf S}^+(0)\right)^{-1}{\bf A}(0)\right]{\bf a}(k).
\end{align*}
Let ${\bf P}(k) \triangleq {\bf I}-\pmb\beta(k)+\pmb\beta(k)\left({\bf S}^+(0)\right)^{-1}{\bf A}(0), {\bf L}^+(0) \triangleq {\bf S}^+(0)-{\bf A}(0)$, and ${\bf a} \triangleq (a_i)_{i=1}^N \neq {\bf 0}$ be a non-trivial weight-balancing solution, i.e., ${\bf a} \in \{\pmb\omega:{\bf L}^+(0)\pmb\omega = {\bf 0}\}$. It follows that ${\bf P}(k)$ and $\bf a$ satisfy the following properties
\begin{enumerate}
\item ${\bf P}(k) \geq {\bf 0}$ since $\beta_i(k) \in [0, 1), \forall i\in \mathcal V, \forall k\in\mathbb Z_+$;
\item ${\bf a} > {\bf 0}$. Since ${\bf L}^+(0)$ is an irreducible singular M-matrix, cf. \cite[C.4.24]{Qu2009}, it follows from \cite[T.4.31]{Qu2009} that 1) ${\bf L}^+(0)$ has rank $N-1$, and 2) $\exists \tilde{\pmb\omega} > {\bf 0}: {\bf L}^+(0)\tilde{\pmb\omega} = {\bf 0}$, and thus ${\bf a} > {\bf 0}$;
\item ${\bf P}(k){\bf a} = {\bf a}, \forall k\in\mathbb Z_+$ since ${\bf S}^+(0){\bf a} = {\bf A}(0){\bf a}$.
\end{enumerate}
Let ${\bf a}_{\rm m} \triangleq (a_{{\rm m},i})_{i=1}^N, {\bf a}_{\rm M} \triangleq (a_{{\rm M},i})_{i=1}^N \in \{\pmb\omega:{\bf L}^+(0)\pmb\omega = {\bf 0}\}$ with $\max_{i\in\mathcal V}a_{{\rm m},i} = \min_{i\in\mathcal V}a_{{\rm M},i} = 1$, ${\bf e}_{\rm m} \triangleq {\bf 1}-{\bf a}_{\rm m} = {\bf a}(0) - {\bf a}_{\rm m} \geq {\bf 0}$ and ${\bf e}_{\rm M} \triangleq {\bf 1}-{\bf a}_{\rm M} \leq {\bf 0}$. It follows that
\begin{align*}
{\bf a}(k+1) 
&{=} {\bf P}(k){\bf a}(k) 
{=} \left[\prod_{t=0}^k{\bf P}(t)\right]\!\!{\bf a}(0) 
{=} \left[\prod_{t=0}^k{\bf P}(t)\right]\!\!({\bf a}_{\rm m}{+}{\bf e}_{\rm m}) \\
&\stackrel{(a)}{=} {\bf a}_{\rm m} + \left[\prod_{t=0}^k{\bf P}(t)\right]{\bf e}_{\rm m}
\stackrel{(b)}{\geq} {\bf a}_{\rm m},
\end{align*}
where $(a)$ follows from ${\bf P}(k){\bf a}_{\rm m} = {\bf a}_{\rm m}, \forall k\in\mathbb Z_+$ and $(b)$ follows from ${\bf P}(k) \geq {\bf 0}, \forall k\in\mathbb Z_+, {\bf e}_{\rm m} \geq {\bf 0}$. Similarly,
\begin{align*}
{\bf a}(k+1) 
= {\bf a}_{\rm M} + \left[\prod_{t=0}^k{\bf P}(t)\right]{\bf e}_{\rm M}
\leq {\bf a}_{\rm M},
\end{align*}
which proves the desired results with
\begin{align}
a_{\min} = \min_{i\in\mathcal V}a_{{\rm m},i}, \quad \text{and} \quad a_{\max} = \max_{i\in\mathcal V}a_{{\rm M},i}. \label{a_globalbound}
\end{align}

\vspace{-0.5cm}
\section{Distributed Quantized Average Consensus} \label{section_consensus}
In this section, we study convergence   of   Algorithm~\ref{alg_dac}. We introduce  the following mild assumptions. 

The first condition is on the number of bits used to quantize the consensus variables at each iteration.  
\begin{assumption} \label{assump_Bc} Let $\left\{B^{\rm (c)}(k)\right\}_{k \in \mathbb{Z}_+}$ be an activation sequence  satisfying $B^{\rm (c)}(k){\in}\{0,1\}$ and $ \sum_{t = nW}^{(n+1)W-1}{B^{\rm (c)}(t)} \geq 1,$ for all $k,n \in \mathbb{Z}_+$ and   some given $W\in \mathbb{Z}_{++}$. 
The number of bits $\big\{B_i^{\rm (c)}(k)\big\}_{k \in \mathbb{Z}_+}$ used by each node $i$ satisfies
\begin{align*}
\begin{cases}
B_i^{\rm (c)}(k)\geq B^{\rm (c)}(k), & \text{if }B^{\rm (c)}(k)=1; \\
B_i^{\rm (c)}(k)= 0, & \text{else.}
\end{cases}
\end{align*}  
\end{assumption}
The above condition is almost the same as the one used in the weight-balancing algorithm {(cf.~A.\ref{assump_Bw})} except for the global upper bound, and can be coupled with it. For example, nodes can communicate for weight-balancing using one bit at odd time slots, and for average consensus using one bit at even time slots, yielding one bit per channel use. Lower effective data rates can be achieved using intermittent communications.

We next introduce the  assumption on $\bar{y}(0)$ and the step-size used in the consensus updates.
\begin{assumption}[Informative $\bar{y}(0)$] \label{assump_sbar}
The average $\bar{y}(0)$ [cf.~(\ref{eq:consensus})] satisfies $\bar{y}(0) \in [q_{\min}, q_{\max}]$.
\end{assumption}
\begin{assumption} \label{assump_alpha_ms} The step-size sequence
   $\{\alpha(k)\}_{k \in \mathbb{Z}_+}$ satisfies:
\begin{align*}
&\alpha(k)>0, \alpha(k+1) \leq \alpha(k), \forall k \in\mathbb{Z}_+, \\
&\sum\limits_{k=1}^{\infty}{\alpha(k)} = \infty, \sum\limits_{k=1}^{\infty}{\alpha(k)^2}< \infty.
\end{align*}
\end{assumption}

 It is important to remark that A.\ref{assump_sbar} neither  requires $y_i(0)$ to be confined within  the quantization range  nor its to be known. This is a major departure from the literature, which calls for  $y_i(0)$ to be within the quantization range -- see, e.g. \cite{Rajagopal2011,Li2011,Wang2011,Thanou2013}. We require instead the   \emph{average} $\bar y(0)$  to fall within the quantization interval $[q_{\min}, q_{\max}]$, which is a less restrictive condition. 
 For example, if nodes are estimating a common unknown parameter $\theta$ via noisy measurements $y_i(0){=}\theta{+}\omega_i$ corrupted by zero mean Gaussian noise $\omega_i$, i.i.d. across nodes, then $\bar{y}(0)$ is the sample mean estimate across the nodes. In this case, a bound on $y_i(0)$ is hard to obtain (theoretically it is unbounded), but the bound of the parameter, $\theta{\in}[\theta_{\min}, \theta_{\max}]$, is known in many cases.
Even worse, $\max_{i\in\mathcal V}|y_i(0)|{\to}\infty$ for $N{\to}\infty$, whereas the sample average $\bar y(0)\to \theta$, so that it becomes more and more informative for large $N$, whereas the initial local measurements become larger and larger. In this example, nodes can simply set $(q_{\min},q_{\max}){=}(\theta_{\min},\theta_{\max})$, so that ${\bar y}(0)$ is informative with high probability.
Herein, we are not interested in non-informative ${\bar y}(0)$, which, as the name suggests, does not provide any information to estimate $\theta$.

We are now ready to state the convergence  of  Algorithm~\ref{alg_dac}.

\begin{theorem} \label{thm_ac}
Let $\left\{\mathbf{y}(k)=(y_i(k))_{i=1}^N\right\}_{k\in \mathbb{Z}_+}$ be the sequence generated by  Algorithm~\ref{alg_dac} under A.\ref{assump_gamma2}-\ref{assump_alpha_ms}. Then: \\
\noindent \texttt{(a) Almost sure convergence}: 
\begin{align}
\mathbb{P} \left(\lim\limits_{k \rightarrow \infty}{{\bf y}(k) = \bar{y}(0)\cdot {\bf 1}} \right) =1. \label{eq_as_main}
\end{align}
\noindent \texttt{(b) Convergence in the moment generating function}:
\begin{align}
\lim\limits_{k \rightarrow \infty}{\mathbb{E}\left[e^{r\Vert{\bf y}(k) - \bar{y}(0) \cdot {\bf 1}\Vert}\right]} = 1, \quad\forall r \in \mathbb R. \label{eq_ms_main}
\end{align}
Furthermore, if $\alpha(k) = \mathcal O(1/k)$ and $\exists m > 0:\alpha(k) \geq m/(k+1), \forall k$, then: \\
\noindent \texttt{(c) Convergence rate}:
\begin{align}
V_k\triangleq\mathbb E[V({\bf y}(k))]
\leq
\begin{cases}
\mathcal O(1/k)&\rho> 1,
\\
\mathcal O(\ln(k)/k)&\rho=1,
\\
\mathcal O(1/k^{\rho})
&\rho<1,
\end{cases}
\end{align} 
where
$\rho{\triangleq}2\xi_1m>0$ with $\xi_1>0$ defined in L.\ref{lemma_ac_dec_term}.
\end{theorem}
\begin{IEEEproof} 
Let $({\bf y}(k),{\bf x}(k),\tilde{\bf y}(k)){\triangleq}\left(y_i(k),x_i(k),\tilde{y}_i(k)\right)_{i=1}^N$ with $\tilde y_i{=} {\rm clip}(y_{i}; q_{\min}, q_{\max})$.
  Using \eqref{update_consensus},
the $y$-updates become
\begin{align}
\!\!\!{\bf y}(k+1)={\bf y}(k)-\alpha(k){\bf L}^{+}(k)\mathbf x(k).
\label{update_y}
\end{align}
and, due to the  {probabilistic} quantization, $\mathbb E[\mathbf x(k)|\mathbf y(k)]{=}\tilde{\mathbf y}(k)$.

To study  the dynamics of the consensus error, we define
\begin{align}
&V({\bf y})\triangleq\Vert{\bf y}-\bar{y}(0){\bf 1}\Vert^2,\label{V}
\end{align} 
and prove that the sequence $\{V({\bf y}_k)\}_{k\in\mathbb Z_+}$ satisfies the conditions of \cite[T.1]{Kar2009}, sufficient to prove our theorem. 

\noindent \textbf{Intermediate results:} We begin by introducing some properties of $V({\bf y})$, instrumental for the sequel of the proof.
\begin{lemma} \label{lemma_ac_error_relation}
In the setting  of T.\ref{thm_ac}, $\forall \mathbf{y}\in \mathbb{R}^N$ there holds
\begin{align*}
&\mathbb E\left[V({\bf y}(k\!+\!1))|{\bf y}(k)={\bf y}\right] \nonumber\\
&= \begin{cases} V({\bf y}) -2\alpha(k){\bf y}^T{\bf L}^+(k)\tilde{\bf y}\\
\quad+\alpha(k)^2{\mathbb E}\!\left[ \Vert {\bf L}^+(k){\bf x}(k)\Vert^2|{\bf y}(k)={\bf y}\right], &\text{if }B^{\rm (c)}(k)=1, \\
V({\bf y}), & \mbox{otherwise}.
\end{cases}
\end{align*}
\end{lemma}
The case $B^{\rm (c)}(k){=}0$ holds trivially since 
${\bf y}(k+1){=}\mathbf y$. Otherwise
($B^{\rm (c)}(k){=}1$) L.\ref{lemma_ac_error_relation} follows
  from the dynamics \eqref{update_y} and the fact that
  $\mathbb E[\mathbf x(k)|\mathbf y(k)]{=}\tilde{\mathbf y}(k)$ with
  probabilistic quantization. To bound these dynamics when $B^{\rm (c)}(k){>}0$, we use the fact that ${\bf y}(k)$ is uniformly bounded within a bounded set $\mathcal S$ with probability 1 (cf. L.\ref{lemma_y_global_bound}) and L.\ref{lemma_ac_dec_term} to obtain
\begin{align}
&\mathbb E\left[V({\bf y}(k+1))|{\bf y}(k)={\bf y}\right] \nonumber \\
&\leq V({\bf y})
-2\alpha(k)\xi_1 V({\bf y})+ 2\alpha(k)\xi_2\Vert {\bf b}(k)\Vert_1 \nonumber\\
&\quad +\alpha(k)^2\Vert {\bf L}^+(k)\Vert_2^2{\mathbb E}\left[ \Vert {\bf x}(k)\Vert^2|\mathbf y(k)={\bf y} \right] \nonumber \\
&\overset{(a)}{\leq} 
V({\bf y})-2\xi_1\alpha(k)\left[V({\bf y})-c(k)\right], \forall {\bf y} \in \mathcal S, \label{err_descent}
\end{align}
where $\xi_1,\xi_2$ are constants defined in L.\ref{lemma_ac_dec_term} and in $(a)$
    we defined \vspace{-0.2cm}
\begin{equation}  
c(k)\triangleq
\frac{\xi_2}{\xi_1}\Vert {\bf b}(k)\Vert_1+\alpha(k)\frac{\xi_3}{2\xi_1}, \label{c_k}
\end{equation}
for some constant $\xi_3{\geq}\Vert {\bf L}^+(k)\Vert_2^2{\mathbb E}\left[ \Vert {\bf x}(k)\Vert^2|\mathbf y(k) = {\bf y} \right]{>}0$. 
Note that the boundedness of ${\bf A}(k)$ (cf. L.\ref{lemma_a_bound}), and thus  of ${\bf L}^+(k)$, and that of  $\mathbf x(k)$   (being the output of a finite rate quantizer), guarantee that $\xi_3{<}\infty$.
 
 We are now ready to prove T.\ref{thm_ac}.

\noindent \textbf{Proof of statement (a):} Define  $\tilde{\mathcal K} \triangleq  \{k :B^{\rm (c)}(k)=1\}$, $\{{\bf y}_{\tilde{\mathcal K}}(k)\} \triangleq \{{\bf y}(\tilde{k})\}_{\tilde{k} \in \tilde{\mathcal K}}, \{{\bf c}_{\tilde{\mathcal K}}(k)\} \triangleq \{{\bf c}(\tilde{k})\}_{\tilde{k} \in \tilde{\mathcal K}}$ and $\{\alpha_{\tilde{\mathcal K}}(k)\} \triangleq \{\alpha(\tilde{k})\}_{\tilde{k} \in \tilde{\mathcal K}}$. It is sufficient to show that $V$ in {\eqref{V}} satisfies the conditions of \cite[T.1]{Kar2009}, namely:
\begin{align*} 
1)\ &\underset{\left\Vert{\bf y}-\bar{y}(0){\bf 1}\right\Vert \geq \epsilon}\inf V({\bf y}) > 0, \forall \epsilon > 0, \\
&V(\bar{y}(0)\cdot{\bf 1}) = 0, \text{ and }
\underset{{\bf y} \rightarrow \bar{y}(0)\cdot{\bf 1}}\limsup V({\bf y}) = 0;
\\
2)\ & \mathbb E\left[V({\bf y}_{\tilde{K}}(k{+}1))|{\bf y}_{\tilde{K}}(k)  {=}    {\bf y}\right]{-}V({\bf y}) \nonumber\medskip \\ & \qquad {\leq} g(k)\left[1{+}V({\bf y})\right]{-}\alpha_{\tilde{\mathcal K}}(k)\phi({\bf y}),
\end{align*}
where $\phi({\bf y})$ is a non-negative function such that
\begin{align*}
\underset{\left\Vert{\bf y}-\bar{y}(0){\bf 1}\right\Vert \geq \epsilon}\inf \phi({\bf y}) > 0,\ \forall\epsilon > 0; 
\end{align*}
  and $\alpha_{\tilde{\mathcal K}}(k)$ and $g(k)$ satisfy 
\begin{align}
&\alpha_{\tilde{\mathcal K}}(k) > 0, \quad 
\sum_{k=1}^\infty{\alpha_{\tilde{\mathcal K}}(k)} = \infty, \label{eq:cond_3a}\\
&g(k) > 0, \quad 
\sum_{k=1}^\infty{g(k)} < \infty.\label{eq:cond_3b}
\end{align}
 
Conditions in 1) are trivially satisfied by definition [cf. (\ref{V})]. To prove the condition in 2), we use ${\bf 1}^\top{\bf y}_{\tilde{K}}(k{+}1) = {\bf 1}^\top{\bf y}_{\tilde{K}}(k)$, L.\ref{lemma_y_global_bound} in App.\ref{append:ac}, and \eqref{err_descent}, 
yielding, 
\begin{align*}
\mathbb E\left[V({\bf y}_{\tilde{K}}(k+1))|{\bf y}_{\tilde{K}}(k) {=}{\bf y}\right]{-}V({\bf y}){\leq} g(k){-}2\xi_1\alpha_{\tilde{K}}(k)V({\bf y}),
\end{align*} 
with $g(k) {=} 2\xi_1\alpha_{\tilde{K}}(k)c_{\tilde{K}}(k)$ and $\phi({\bf y}) {=} 2\xi_1V({\bf y})$. Moreover,
\begin{align*}
\sum_{k \geq 0}{g(k)}\leq 2\xi_1\sqrt{\left[\sum_{k \in \tilde{\mathcal K}}\alpha(k)^2\right]\left[\sum_{k \in \tilde{\mathcal K}} c(k)^2\right]}\overset{(a)}{<}\infty,
\end{align*}  
where $(a)$ we used $\sum_{k \in \tilde{\mathcal K}}{\alpha(k)^2}{<}\infty$ (cf. A.\ref{assump_alpha_ms}); and $\sum_{k \in \tilde{\mathcal K}}{c(k)^2}{<}\infty$, due to  \eqref{c_k}, A.\ref{assump_alpha_ms}, and T.\ref{thm_wb}. Therefore, the condition in 2) holds. 
 
 Overall,  we have shown that all the conditions of \cite[T.1]{Kar2009} are satisfied, implying   $$\mathbb{P} \Bigr(\underset{k \rightarrow \infty, k{\in}\tilde{\mathcal K}}\lim{{\bf y}(k)=\bar{y}(0)\cdot {\bf 1}} \Bigr)=1.$$   Since $|\tilde{\mathcal K}|=\infty$ and ${\bf y}(k+1){=}{\bf y}(k),$ for all $k \notin \tilde{\mathcal K}$,  statement (a) of the theorem   follows.

\noindent \textbf{Proof of statement (b):} 
Since $\Vert{\bf y}(k)-\bar{y}(0){\bf 1}\Vert^r{<}\infty$ and $\mathbb{P} \Bigr(\underset{k \rightarrow \infty}\lim{\Vert{\bf y}(k)-\bar{y}(0){\bf 1}\Vert^r = 0} \Bigr)=1$, for all   $r\in\mathbb{Z}_{++}$ (recall that $\left\vert y_i(k){-}\bar{y}(0)\right\vert$ is bounded for all $i{\in}\mathcal{V}$, cf. L.\ref{lemma_y_global_bound} in App.\ref{append:ac}), it follows from the dominated convergence theorem (cf. \cite[T.1.6.7]{Durrett2013}) that $\lim_{k\to \infty}\mathbb E\left[\Vert{\bf y}(k)-\bar{y}(0){\bf 1}\Vert^r\right] = 0, \forall r >0$, which implies statement (b). 

\noindent \textbf{Proof of statement (c):}
For simplicity, we assume that $B^{\rm (c)}(k) = 1, \forall k \in \mathbb Z_+$, and the proof can be easily generalized to the case that $B^{\rm (c)}(k)$ satisfying Assumption \ref{assump_Bc}. Since $\alpha(k){=}\mathcal O(1/k)$ and $\Vert\mathbf b(k)\Vert_1{=}\mathcal O(1/k)$,
it follows that $c(k){=}\mathcal O(1/k)$, and there exist
 $M>0,C>0$ such that
$$\alpha(k)\leq M/(k+1),\ c(k)\leq C/(k+1),\forall k.$$
Under the conditions of the theorem, \eqref{err_descent} holds, which implies
\begin{align}
&V_{k+1}
\leq
\Big(1-\frac{\rho}{k+1}\Big)V_k+\gamma\frac{1}{(k+1)^2},\ \forall k,
\end{align}
where $\gamma{\triangleq}2\xi_1MC$.
Let $\bar k{\triangleq}\lceil\rho\rceil{-}1$.
By induction, we can show
\begin{align}
\nonumber
&V_{k}
\leq
\left[(\bar k+1)V_{\bar k}+\gamma\right]
\beta(\bar k,k)
+\gamma\sum_{t=\bar k+1}^{k-1}\beta(t,k),\ \forall k>\bar k,
\\&
\text{where }\beta(t,k)
\triangleq
\frac{1}{(t+1)^2}
\prod_{i=t+1}^{k-1}\Big(1-\frac{\rho}{i+1}\Big).
\end{align}
Note that
$$
\ln\beta(t,k)
\leq
\int_{t+2}^{k+1}
\ln\left(1-\frac{\rho}{x}\right)\mathrm dx
-2\ln(t+1)
$$$$
=
(k+1-\rho)\ln(1-\rho/(k+1))
-(t+2-\rho)\ln(1-\rho/(t+2))
$$$$
+2\ln((t+2)/(t+1))
+(\rho-2)\ln(t+2)
-\rho\ln(k+1).
$$
The first three terms are bounded,
  since
  $\ln((n{+}2){/}(n{+}1)){\to}0$
  and
  $(n{-}\rho)\ln(1{-}\rho/n){\to}{-}\rho$
   for $n{\to}\infty$. It follows that
$
\beta(t,k)
{\leq}
e^{Q}(t{+}2)^{\rho-2}(k{+}1)^{-\rho}
$,
for some $Q<\infty$.
 Letting $S_k\triangleq(k+1)^{-\rho}\sum_{t=\bar k+3}^{k+1}
t^{\rho-2}$,
it follows that
\begin{align}
V_{k}
\leq
A(k+1)^{-\rho}
+\gamma e^{Q}S_k,\ \forall k>\bar k,
\end{align}
for some $A{<}\infty$.
To conclude, note that $A(k+1)^{-\rho}{=}\mathcal O(k^{-\rho})$,
$$
S_k
\leq
(k+1)^{-\rho}
\int_{1}^{k+2}x^{\rho-2}\mathrm dx
=
\begin{cases}
\mathcal O(1/k),
&\rho>1,\\
\mathcal O(\ln(k)/k),
&\rho=1,\\
\mathcal O(1/k^\rho),
&0<\rho<1,
\end{cases}
$$
which proves the desired result.

\end{IEEEproof}

\section{Numerical Results} \label{section_simulation}
\begin{figure}[t] 
	\centering
	\includegraphics[width = 0.5\linewidth]{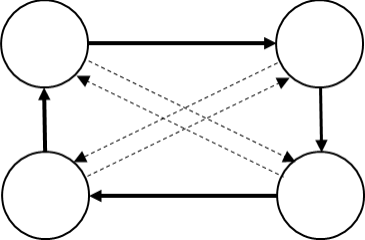} 
	\caption{Illustration of the random graph model for $N=4$, where dashed arrows represent potential directed links depending on the realizations.}	 \label{fig_graph_sim}	
\end{figure}

In this section, we present some numerical results to validate our theoretical findings on strongly connected digraphs with $N = 50$ nodes constructed by the following procedure: a directed ring links all the nodes, to ensure strong connectivity (cf. Fig. \ref{fig_graph_sim}). Then directed edges are randomly added, with probability $0.2$ on each pair of nodes.\vspace{-0.2cm}

\subsection{Quantized weight-balancing}
We adopt \eqref{gamma_ex} for $\{\gamma(k)\}_{k \in \mathbb{Z}_+}$. We compare the total imbalance $\Vert {\bf b}(k)\Vert_1$ of our proposed scheme with the integer weight-balancing and real weight-balancing schemes in \cite{Rikos2014}.
The real weight-balancing scheme uses real valued communications; the integer weight-balancing scheme uses \emph{unicast} transmissions to each of its out-neighbors to communicate the associated edge weight, and cannot use a prescribed number of bits. As we will see numerically, these features allow the scheme to converge to a weight-balanced solution within finite time.
In contrast, our scheme uses \emph{broadcast} communications with a prescribed number of bits per channel use, which in general does not guarantee convergence within finite time. The simulation results are averaged over 100 graph realizations.

\begin{figure}[t] 
	\centering
	\includegraphics[width = \linewidth]{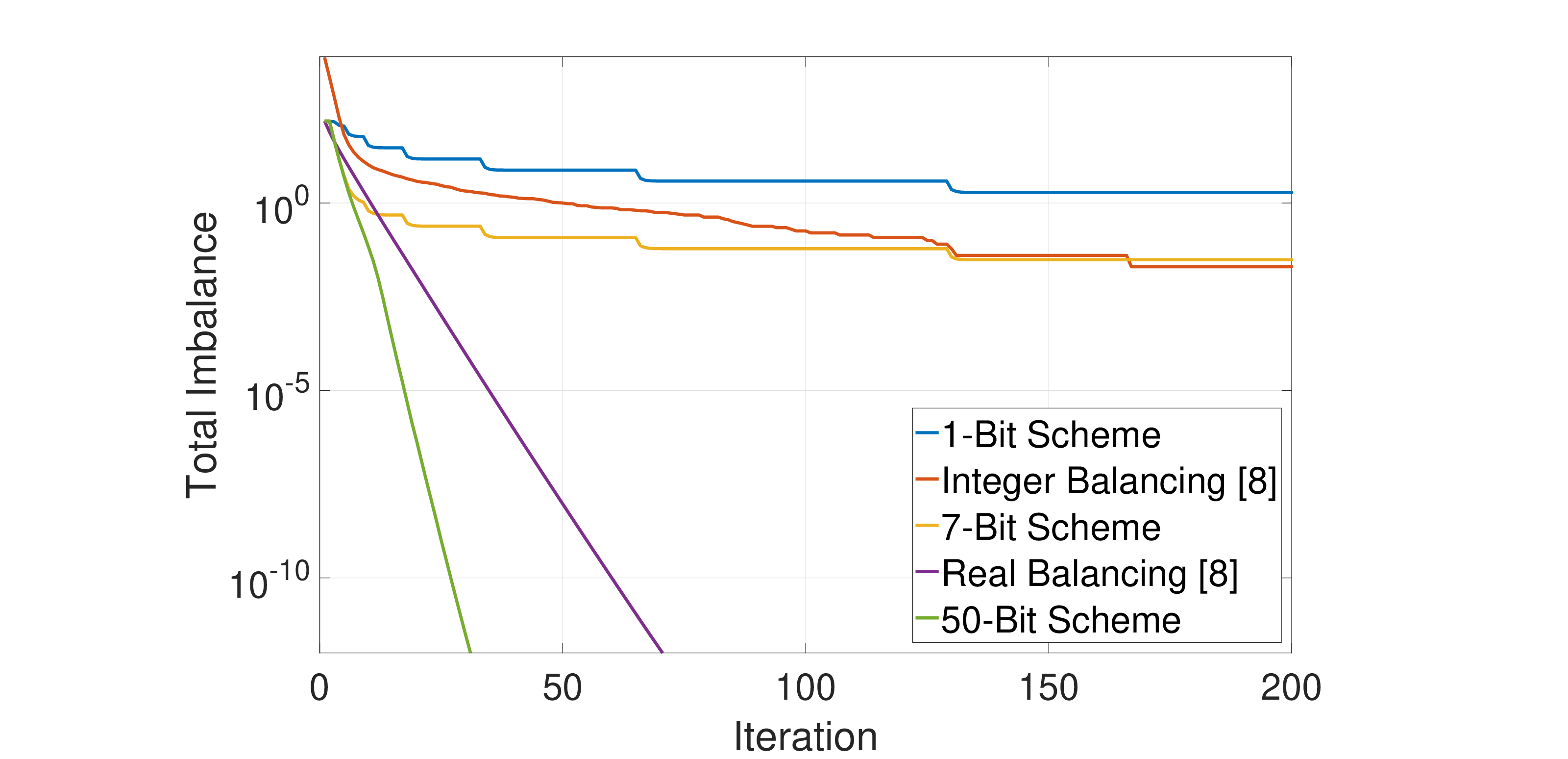} 
	\caption{Quantized weight-balancing problem: Total imbalance $\Vert {\bf b}(k)\Vert_1$ the propsoed algorithm with 1-bit , 7-bit, and 50-bit, as well as the integer and real weight-balancing schemes \cite{Rikos2014}.\vspace{-0.5cm}} \label{fig_imb}	
\end{figure}

Fig. \ref{fig_imb} shows the total imbalance of Algorithm \ref{alg_dwb} with 1-bit and 5-bit of information exchange, as well as of the other two benchmark schemes. Note that in the integer weight-balancing scheme, the maximum weights in the 100 realizations are between 88 to 250, implying that 7 to 8 bits are required \emph{per edge} per timeslot, which implies $7$ or $8\times (1+0.2\times 49) = 75.6$ or $86.4$ bits per node per timeslot. It is shown that the metric $\Vert {\bf b}(k) \Vert_1$ is non-increasing for the proposed schemes, which is consistent with our analytical results (cf. L.\ref{lemma_decr_event}). In addition, one can see that the curve of $\Vert {\bf b}(k) \Vert_1$ can be partitioned into nearly flat and steep line segments, for both schemes. The rationale behind this behavior is that the total imbalance decreases only when decreasing events occur (steep line segments); in between, the imbalance may be transferred within the network, but without causing the total imbalance to decrease. Compared with the two benchmark schemes, it shows that the proposed scheme with 50 bits outperforms the real weight-balancing scheme \cite{Rikos2014}, which requires infinite rate communications. On the other hand, The comparison between the proposed 7-bit scheme and the integer weight-balancing scheme shows that, initially, the proposed scheme has better performance. However, as noted earlier, the integer weight-balancing scheme later outperforms the proposed 7-bit scheme since it is guaranteed to converge to a weight-balanced solution in finite timeslots.\vspace{-0.3cm}
\subsection{Quantized average consensus}
We compare our proposed algorithm with the following state-of-the-art schemes:
1) Q-Push-Sum, where we straightforwardly apply the finite-bit probabilistic quantization to the original push-sum algorithm in \cite{Kempe2003}, i.e., $z_i(0){=}s_i, \psi_i(0){=}1, \forall i{\in}\mathcal V$,
\begin{align*}
\psi_i(k+1) &= \psi_i(k)+\alpha(k)\!\!\sum\limits_{j \in \mathcal{N}_i^-}{\!\!\! a_{ij}(k)\,\left[\mathcal Q \left(\psi_j(k)\right)-\mathcal Q \left(\psi_i(k)\right)\right]}; \\
z_i(k+1) &= z_i(k)+\alpha(k)\!\!\sum\limits_{j \in \mathcal{N}_i^-}{\!\!\! a_{ij}(k)\,\left[\mathcal Q \left(z_j(k)\right)-\mathcal Q \left(z_i(k)\right)\right]};
\end{align*}
and $y_i(k) = \frac{z_i(k)}{\psi_i(k)}$ is the estimate of the initial average,
where $\mathcal Q \left(\bullet\right)$ is the quantization defined in \eqref{quantizer}; note that Q-Push-Sum can be regarded as the generalization of \cite{Rikos2018} to real valued initialization and finite rate communications; 2) Q-Run-Avg, where we apply the finite-bit probabilistic quantization to the algorithm in \cite{Zhu2015};\footnote{Note: this algorithm requires $\mathcal O(N)$ memory space to store the estimate of the eigenvector of graph Laplacian at each node.} 3) Q-Monte-Carlo, where we apply the $B$-bit quantization
\begin{align*}
r_\beta(x; B) = (1+\beta)^{\max\left\{\lfloor \log_{1+\beta}x\rfloor, 2^B\right\}}
\end{align*}
to the Monte-Carlo based algorithm in \cite[Section 4]{Bernadette2018}. In this algorithm nodes exchange quantized random values sampled from the exponential distribution with parameter related to their current states. Note that exact convergence can be achieved by \cite{Zhu2015,Bernadette2018} using infinite-bit quantized communications, \cite{Kempe2003} using real value communications, and \cite{Rikos2018} using integer value communications.
However, there is no theoretical guarantee for all these benchmark schemes with finite bit quantization: our proposed scheme is the first algorithm solving the distributed average consensus over unbalanced digraphs with a prescribed finite rate communications.

%
We adopt the mean square error (MSE) 
${\rm MSE}(k)=V({\bf y}(k))/N$
 as defined in \eqref{V} as performance metric.
The simulation results are averaged over 100 graph realizations and 100 initial value realizations, i.e., totally 10000 realizations.

For the proposed algorithm, we adopt: $q_{\min}{=}0, q_{\max}{=}1$, $B_i^{\rm (w)}(k){=} B_i^{\rm (c)}(k){=}50, \forall i, \forall k$; \eqref{gamma_ex} is adopted for $\{\gamma(k)\}_{k \in \mathbb{Z}_+}$, and $\alpha(k){=}(k+1)^{-1}, \forall k{\in}\mathbb{Z}_+$, which satisfies A.\ref{assump_alpha_ms};
 for Q-Push-Sum we use 50 bits and $q_{\min}{=}0, q_{\max}{=}N$ to quantize $\psi$, and 50 bits and $q_{\min}{=}0, q_{\max}{=}1$ to quantize $z$; for Q-Run-Avg we quantize each element of ${\bf z} \in {\mathbb R}^N$ (the estimate of the left eigenvector at 0 of graph Laplacian constructed by a row stochastic weight matrix) using 4 bits (i.e., totally $4{\times}N{=}80$ bits required for quantizing $\bf z$) and $q_{\min}{=}0, q_{\max}{=}N^\kappa$ with $\kappa{=}1.15$, and $y$ is quantized using 20 bits and $q_{\min}{=}0, q_{\max}{=}1$; for Q-Monte-Carlo, we use 50 bits to quantize both $X$ and $Y$, and other parameters are: $a{=}0, b{=}1, \varepsilon{=}10^{-3}$. Note that the communication resource budget per node per timeslot is 100 bits in all schemes.

\begin{figure}[t] 
	\centering
	\includegraphics[width = \linewidth]{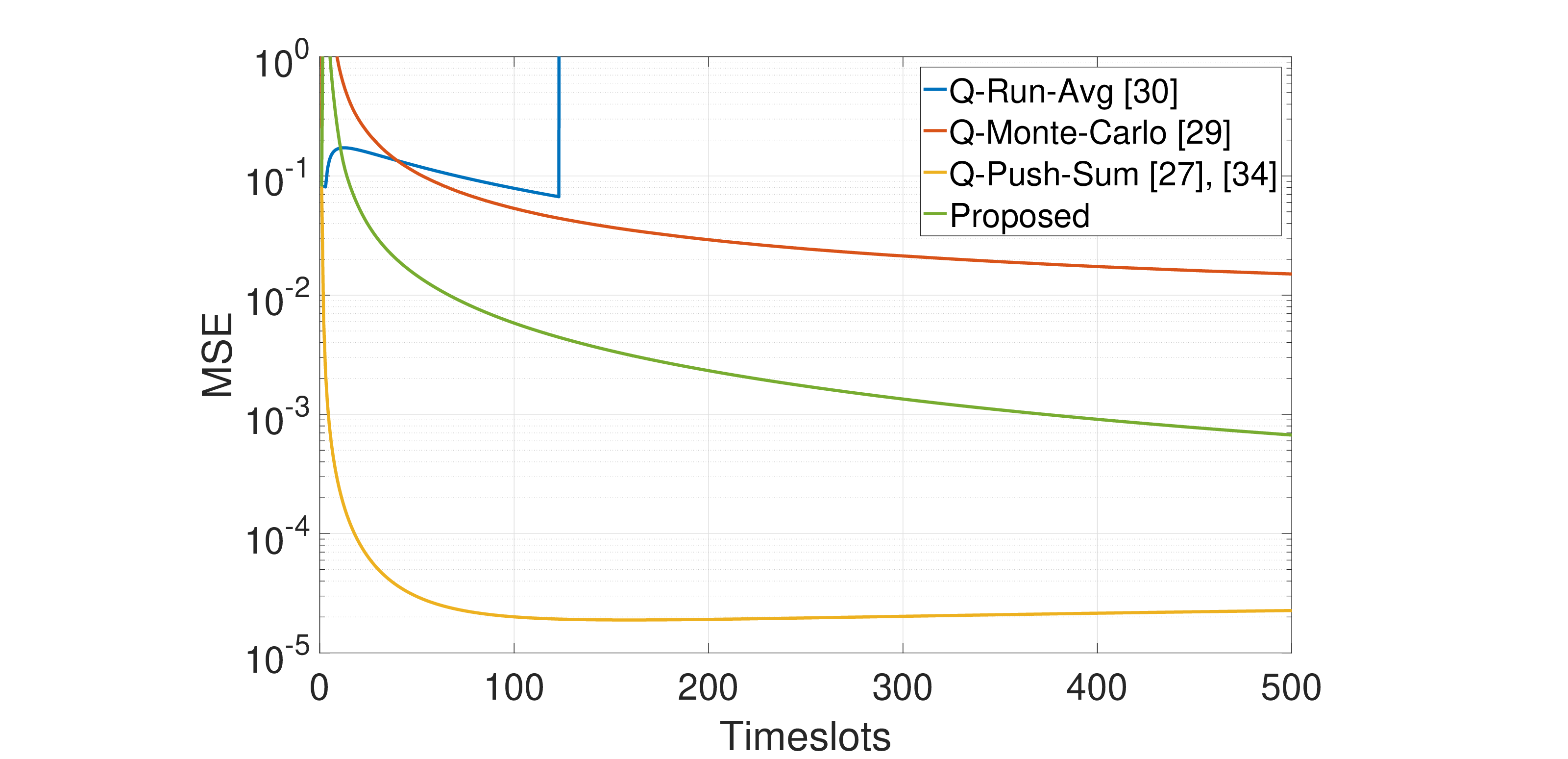}
	\caption{Quantized average consensus problem: MSE of average consensus algorithms average over 10000 realizations.\vspace{-0.5cm}} 		\label{fig_mse}
\end{figure}

\begin{figure}[t] 
	\centering
	\includegraphics[width = \linewidth]{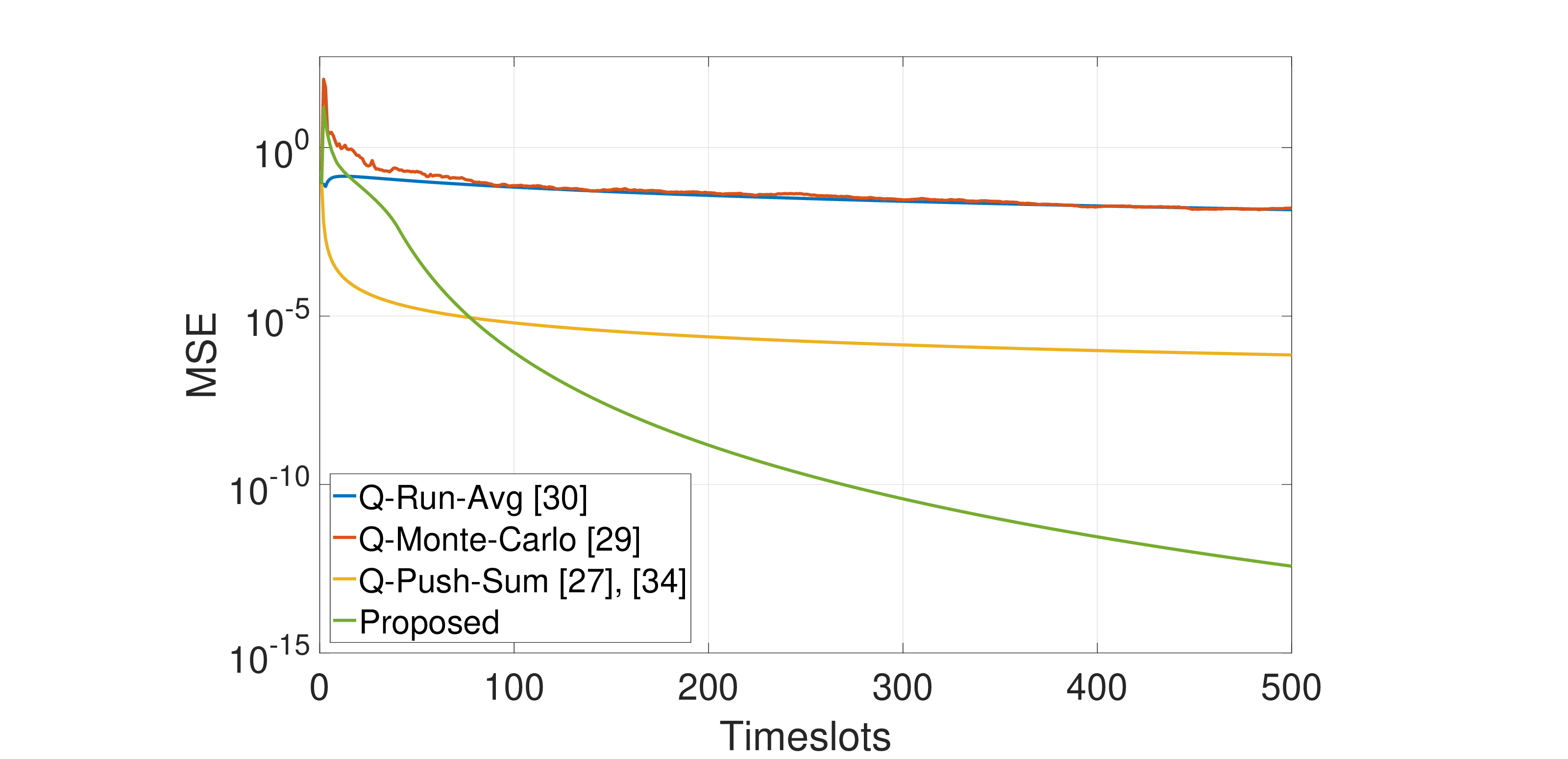}
	\caption{Quantized average consensus problem: MSE of average consensus algorithms for a particular realization.\vspace{-0.5cm}}		\label{fig_mse2} \vspace{0.45cm}
\end{figure}

Fig. \ref{fig_mse} shows the MSE performance of Algorithm \ref{alg_dac} as well as other benchmark schemes. It is shown that only the proposed scheme and the Q-Monte-Carlo are reaching the average consensus, among all finite rate schemes. Note that Q-Run-Avg and Q-Push-Sum seem also converge for some realizations, cf. Fig. \ref{fig_mse2}. However, only the proposed scheme has theoretical convergence guarantees.

\begin{figure}[t] 
	\centering
	\includegraphics[width = \linewidth]{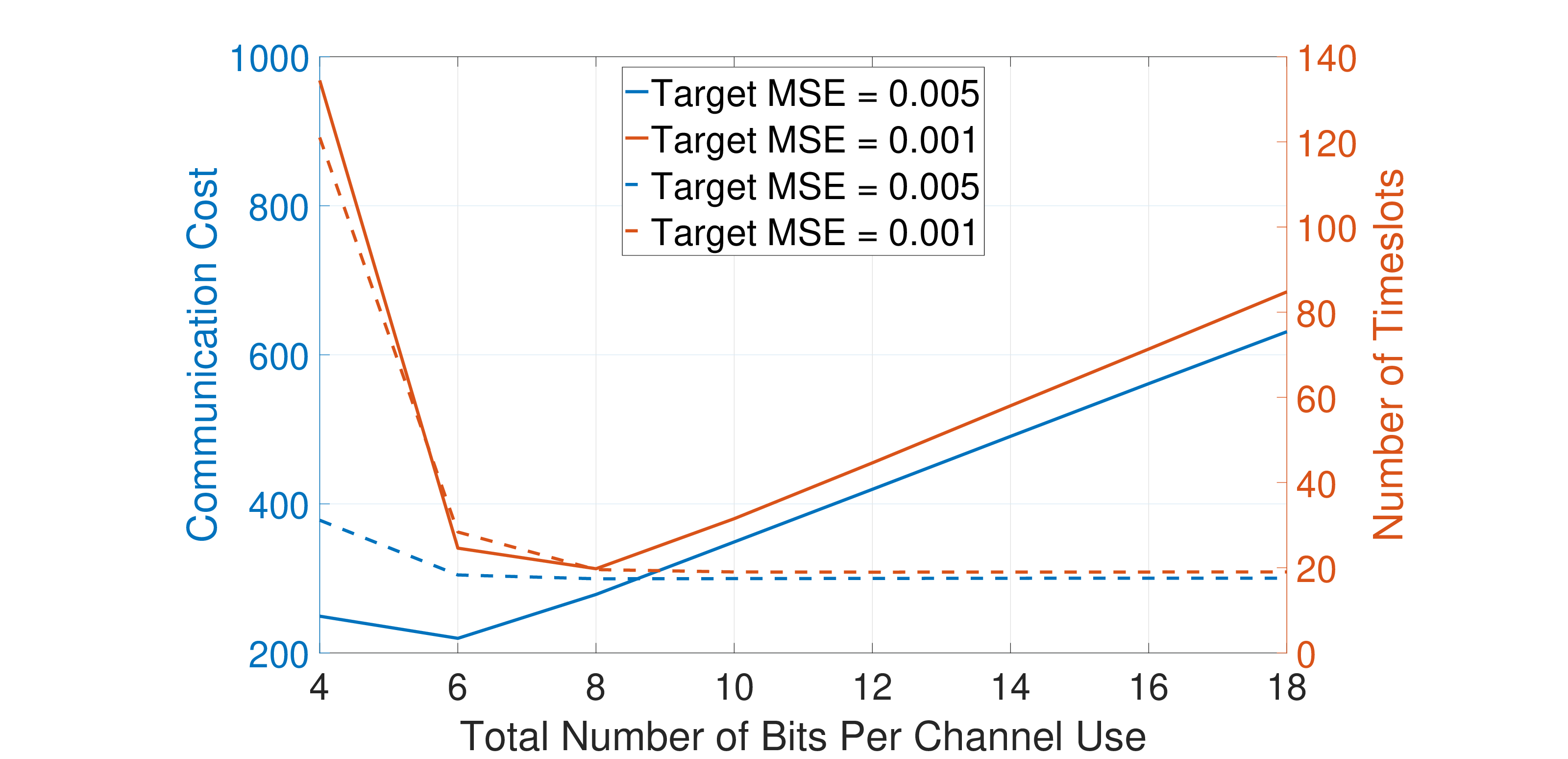}
	\caption{Quantized consensus problem: Communication cost (left y-axis, solid lines) and number of timeslots (right y-axis, dashed lines) needed to reach the target MSE, versus the total number of bits per channel use.\vspace{-0.5cm}}		\label{fig_commcost}
\end{figure}

Fig. \ref{fig_commcost} shows the communication cost (left y-axis) and delay (right y-axis) needed by Algorithm \ref{alg_dac} to reach a target  MSE of $1\times 10^{-3}$ and $5\times 10^{-3}$, versus the total number of bits per channel use.  The communication cost is defined as the product of the total number of bits per node per timeslot and the number of timeslots. For each parameter setting, we run 50 graph realizations and 10 initial value realizations. To avoid the average results affected by the outliers, we select the best $95\%$ of results to perform averaging.
We observe that increasing the total number of bits reduces the number of timeslots required. On the other hand, there exists an optimal number of bits that minimizes the communication cost. Using more bits does not appear to be beneficial, since the communication cost becomes larger, and it is only marginally compensated by the reduction of the number of timeslots required.\vspace{-0.4cm}  

\section{Conclusions} \label{section_conclusion}
In this paper, we introduced a novel distributed algorithm that solves the weight-balancing problem using only quantized information and simplex communications. Building on this scheme, a second contribution of the paper was a novel distributed average consensus algorithm over non-balanced digraphs that uses quantized simplex communications. Convergence of the algorithm was proved using a novel line of analysis, based on  a novel metric inspired by the positional  system representation and a  new step-size rule. Finally, numerical results validated our theoretical findings.

\appendices
\def\thesubsectiondis{\Roman{subsection}.} 

\section{Intermediate Results in the  Proof of Theorem~\ref{thm_wb}} \label{append:wb}
\subsection{Preliminary definitions and results} \label{append:notation}
Throughout the proof, we write the updates of $S_i^+(k), S_i^-(k)$, $b_i(k)$ of Algorithm~\ref{alg_dwb} as
\begin{align}
S_i^+(k+1) &\triangleq \sum\limits_{j = 1}^N{a_{ji}(k+1)} \overset{\eqref{a_ij}}{=} \sum\limits_{j \in \mathcal{N}_i^+}{\left[a_{ji}(k)+\gamma(k)n_i(k)\right]} \nonumber\\
&= S_i^+(k)+d_i^+\gamma(k)n_i(k), \label{S_i+} \\
S_i^-(k+1) &\triangleq \sum\limits_{j = 1}^N{a_{ij}(k+1)} \overset{\eqref{a_ij}}{=} \sum\limits_{j \in \mathcal{N}_i^-}{\left[a_{ij}(k)+\gamma(k)n_j(k)\right]} \nonumber\\
&
= S_i^-(k)+\sum\limits_{j \in \mathcal{N}_i^-}{\gamma(k)n_j(k)}, \label{S_i-} \\
b_i(k+1) &\triangleq S_i^-(k+1)-S_i^+(k+1) \nonumber\\
&= b_i(k)-d_i^+\gamma(k)n_i(k)+\sum_{j\in \mathcal{N}_i^-}\!\!\!\gamma(k)n_j(k). \label{b_i}
\end{align}

\begin{lemma} \label{lemma_node_flow}
Given 
$\mathcal{V}^+(k)$ and $\mathcal{V}^-(k)$, defined in  \eqref{V+V-}, it holds: 
 $$\mathcal{I}(\mathcal D_k){=}0\,\,\Rightarrow\,\,\mathcal{V}^+(k+1) {\supseteq} \mathcal{V}^+(k) \text{ and }   \mathcal{V}^{-}(k+1) {\supseteq} \mathcal{V}^{-}(k).$$ 
\end{lemma}
\begin{IEEEproof}
Let $\mathcal I\{\mathcal D_k\}{=}0$ and consider $i{\in}\mathcal V^-(k)$. Then, \eqref{signal_wb_alg_ac}  and $b_i(k){<}0$ imply $\gamma(k)n_i(k){>}b_i(k)/d_i^+$;
$\mathcal I\{\mathcal D_k\}{=}0$ (see \eqref{event_D})
implies $n_{j}(k){\leq}0, \forall j{\in}\mathcal{N}_i^-$. 
$b_i(k{+}1){<}0$ then follows from \eqref{b_i}, so that
$i{\in}\mathcal V^-(k{+}1)$; hence  $\mathcal{V}^-(k{+}1){\supseteq}\mathcal{V}^-(k)$.
$\mathcal{V}^+(k{+}1){\supseteq}\mathcal{V}^+(k)$ follows
from a similar argument on $i{\in}\mathcal V^+(k)$.
\end{IEEEproof}

\begin{lemma} \label{lemma_T-}
 $\forall {k} \in \mathbb{Z}_+,\exists\tau \geq  k:\Vert {\bf b}(\tau)\Vert_{1} < N^2\gamma(\tau)$.
\end{lemma}
\begin{IEEEproof}
We prove it by contradiction. Let $\bar{T}_0 \triangleq (2W-1)\bar{U}$. Suppose $\exists k \in \mathbb{Z}_+:\Vert {\bf b}(\tau)\Vert_1 \geq N^2\gamma(\tau)$, for all $\tau\geq k$. Invoking  Corollary~\ref{corollary_minimal_decrease}.{\it(ii)} recursively $m$ times and taking $m\to\infty$ yields 
\begin{align*}
0\leq\left\Vert {\bf b}\big(k\!+\!m\,\bar{T}_0\big)\right\Vert_1
\!\!\leq \!
\left\Vert {\bf b}\big(k\big)\right\Vert_1
\!\!-\!2\sum_{n=1}^{m}\!\gamma\big(k\!+\!n\,\bar{T}_0\big),
\end{align*}
a contradiction since $\sum_{n=1}^{m}\gamma\left(k+n\bar T_0\right)\to\infty$ due to \eqref{eq:lower_upper_gamma}.
   \hfill 
\end{IEEEproof}

\subsection{Proof of Lemma \ref{lemma_a_bound}}\label{sec_proof_lemma_A}
By the definition of Cauchy sequence applied to each entry of ${\bf A}(k)$, we need to prove that, $\forall \epsilon>0,\exists k_\epsilon\in\mathbb Z_+$ such that $\max_{i,j}|a_{ij}(m)-a_{ij}(n)|<\epsilon,\forall m,n\geq k_\epsilon$.
To this end, let $\epsilon{>}0$ and define $k_\epsilon$ as\footnote{Note that $k_\epsilon^{(i)}<\infty$ since $\Vert {\bf b}(k) \Vert_1\to 0$ and $\gamma(k)\to 0$, see \eqref{eq:lower_upper_gamma}.} 
\begin{align}
\label{keps}
k_\epsilon
 = \min \{k: \Vert {\bf b}(k) \Vert_1 < \epsilon/(2\bar{U}),\epsilon/(4\bar{U})\}.
\end{align}
Since $a_{ji}(k)$ is updated only at update or decreasing events,
using \eqref{a_ij} recursively, we infer  
$$
|a_{ji}(k)-a_{ji}(k_\epsilon)|
\leq \sum_{\ell=k_\epsilon}^{\infty}\mathcal I\{\mathcal U_\ell\vee\mathcal D_\ell\}\gamma(\ell),\ \forall k\geq k_\epsilon.
$$
With $t_l$ defined as in \eqref{t_l} and letting $L_\epsilon\triangleq \min\{l\in\mathbb Z_+:t_l\geq k_\epsilon\}$, we can further upper bound
\begin{align}
|a_{ji}(k)-a_{ji}(k_\epsilon)|
\leq &
\gamma(k_\epsilon)\sum_{\ell=k_\epsilon}^{t_{L_\epsilon}-1}\mathcal I\{\mathcal U_\ell\vee\mathcal D_\ell\}
\\&
+\sum_{l=L_\epsilon}^\infty\gamma(t_l)\sum_{\ell=t_l}^{t_{l+1}-1}\mathcal I\{\mathcal D_\ell\vee\mathcal U_\ell\}.
\nonumber
\end{align}
Since there are at most $\bar{U}$ update events between the two consecutive decreasing events 
at times $t_l$ and $t_{l+1}$ (cf. C.\ref{corollary_minimal_decrease}),
it follows that $\sum_{\ell=t_l}^{t_{l+1}-1}\mathcal I\{\mathcal D_\ell\vee\mathcal U_\ell\}\leq \bar{U}$, hence
\begin{align}
\label{zzzz}
|a_{ji}(k)-a_{ji}(k_\epsilon)|
\leq 
\left[\gamma(k_\epsilon)
+\sum_{l=L_\epsilon}^\infty\gamma(t_l)\right]\bar U.
\end{align}
To bound $\sum_{l=L_\epsilon}^\infty\gamma(t_l)$, we
apply recursively L.\ref{lemma_decr_event},
\begin{align}
\nonumber
-\Vert {\bf b}(k_\epsilon)\Vert_1=
\sum_{\ell=k_\epsilon}^\infty
[\Vert {\bf b}(\ell+1)-\Vert {\bf b}(\ell)\Vert_1]
\\
\leq
-2\sum_{\ell=k_\epsilon}^\infty\gamma(\ell)\mathcal I\{\mathcal D_\ell\}
=
-2\sum_{l=L_\epsilon}^\infty\gamma(t_l),
\label{dfghs}
\end{align}
hence $\sum_{l=L_\epsilon}^\infty\gamma(t_l)\leq\Vert {\bf b}(k_\epsilon)\Vert_1/2
\leq \epsilon/(2\bar{U})$.
By combining \eqref{zzzz} with \eqref{dfghs} and \eqref{keps}, we finally obtain, $\forall k\geq k_\epsilon$,
\begin{align*}
|a_{ji}(k)-a_{ji}(k_\epsilon)|
\leq 
\left[\gamma(k_\epsilon)
+\Vert {\bf b}(k_\epsilon)\Vert_1/2\right]\bar U
<\epsilon/2,
\end{align*}
and, $\forall n,m\geq k_\epsilon$,
\begin{align*}
|a_{ji}(n){-}a_{ji}(m)|\leq|a_{ji}(n){-}a_{ji}(k_\epsilon)|
+|a_{ji}(m)-a_{ji}(k_\epsilon)|{<}\epsilon,
\end{align*}
which proves that $\{{\bf A}(k)\}_{k \in \mathbb{Z}_+}$ is a Cauchy sequence. \vspace{-0.2cm}

\subsection{Proof of Lemma \ref{lemma_decr_event}} \label{append:lemma_decr_event}
We first introduce the following intermediate result.
\begin{lemma} \label{lemma_eps_gamma}
Let $\{{\bf b}(k)\}_{k \in \mathbb{Z}_+}$ be the sequence generated by Algorithm \ref{alg_dac}. Then, $b_i(k)/\bar\gamma(k){\in}\mathbb{Z},\forall i \in \mathcal{V}$ and $k \in \mathbb{Z}_+$.\vspace{-0.2cm}
\end{lemma}
\begin{IEEEproof}
We prove this lemma by induction using \eqref{b_i}. 
The induction hypothesis holds at $k{=}0$, since $a_{ij}(0){=}1,\forall (j,i)\in \mathcal{E}$ and $\bar\gamma(0){=}2^{1-B_{\max}}$ (cf. A.\ref{assump_gamma2}).
 Suppose that it holds at $k{\geq}0$, i.e., $b_i(k)/\bar\gamma(k){\in}\mathbb{Z},\forall i$.
Then, since $\bar\gamma(k+1){=}\bar\gamma(k)/m(k)$, with $m(k){\in}\{1,c_1\}{\subset} \mathbb{Z}_{++}$ (A.\ref{assump_gamma2}), it follows that $b_i(k)/\bar\gamma(k+1)=m(k)b_i(k)/\bar\gamma(k){\in}\mathbb{Z}$ and $\gamma(k)n_i(k)/\bar\gamma(k+1)=m(k)n_i(k)2^{B_{\max}-1}{\in}\mathbb Z$. Therefore, by \eqref{b_i}, one can infer that $b_i(k+1)/\bar\gamma(k+1)\in\mathbb Z$, proving the induction step and completing the proof.
\end{IEEEproof}

\noindent \textbf{Proof of Lemma \ref{lemma_decr_event}:}
Let $i \in \mathcal{V}$ and $k \in \mathbb{Z}_+$. Using \eqref{b_i}, we find
\begin{align}
\label{xzdfg}
\frac{|b _i(k+1)|}{\gamma(k)}=\left\vert \frac{b_i(k)}{\gamma(k)} -d_i^+n_{i}(k)+\sum_{j\in \mathcal{N}_i^-}n_{j}(k)\right\vert.
\end{align}
It will be useful to note that, as can be seen from \eqref{signal_wb_alg_ac},
$b_i(k)$, $n_i(k)$ (possibly, $n_i(k)=0$) and $b_i(k)-d_i^+n_i(k)$ have the same signs, yielding the following inequality for all $i$,
\begin{align}
\label{qwer}
\frac{|b _i(k+1)|}{\gamma(k)}\stackrel{\triangle}{\leq} \frac{|b_i(k)|}{\gamma(k)}-d_i^+|n_{i}(k)|+\sum_{j\in \mathcal{N}_i^-}\vert n_{j}(k)\vert,
\end{align}
where $\triangle$ stands for the triangle inequality. We now distinguish the two cases  $\mathcal I\{\mathcal D_k\}{=}0$ and $\mathcal I\{\mathcal D_k\}{=}1$.
If $\mathcal I\{\mathcal D_k\} = 0$,
from the negation of $\mathcal D_k$ in \eqref{event_D} and from \eqref{signal_wb_alg_ac}, it follows that
$b_i(k), n_i(k), b_i(k)-d_i^+n_i(k)$, and $n_j(k), \forall j \in \mathcal N_i^-$ have the same signs, so that \eqref{qwer} holds with equality, $\forall i$.

Conversely, if $\mathcal I\{\mathcal D_k\}{=}1$, there exists $\ell$
and $(j_1,j_2){\in}[\mathcal{N}_\ell^-]^2$
 such that either 1) $b_{\ell}(k){\neq}0$, $n_{\ell}(k){\neq}0$ and 
 ${\rm sgn}(n_{j_1}(k)){=}-{\rm sgn}(b_{\ell}(k))$;
 or 2) $b_{\ell}(k){=}0$, $n_{\ell}(k){=}0$ and ${\rm sgn}(n_{j_1}(k)){=}-{\rm sgn}(n_{j_2}(k))$.
 In the first case ($b_{\ell}(k)\neq 0$, $n_{\ell}(k)\neq 0$ and 
 ${\rm sgn}(n_{j_1}(k))=-{\rm sgn}(b_{\ell}(k))$), we bound \eqref{xzdfg} as
 \begin{align*}
 \nonumber
&\frac{|b _{\ell}(k+1)|}{\gamma(k)}\stackrel{\triangle}{\leq}\left\vert \frac{b_{\ell}(k)}{\gamma(k)} -d_{\ell}^+n_{\ell}(k)+n_{j_1}(k)\right\vert
+\!\!\!\!\!\!\!\!\!\sum_{j\in \mathcal{N}_\ell^-,j\neq j_1}\!\!\!\!\!\vert n_{j}(k)\vert
\\
&
{\leq}\left\vert \frac{|b_\ell(k)|}{\gamma(k)} {-}d_\ell^{+}|n_{\ell}(k)|-|n_{j_1}(k)|\right\vert
+\sum_{j\in \mathcal{N}_\ell^-,j\neq j_1}\vert n_{j}(k)\vert
\\
\nonumber
&{=} 
\frac{|b_\ell(k)|}{\gamma(k)}{-}d_\ell^+|n_{\ell}(k)|
-2\min\left\{|n_{j_1}(k)|,\frac{|b_\ell(k)|}{\gamma(k)}{-}d_\ell^+|n_{\ell}(k)|\right\}\\&
\quad+\sum_{j\in \mathcal{N}_\ell^-}\vert n_{j}(k)\vert.
\nonumber
\end{align*}
 In the second case ($b_{\ell}(k){=}0$, $n_{\ell}(k){=}0$, 
 ${\rm sgn}(n_{j_1}(k)){=}-{\rm sgn}(n_{j_2}(k))$ and, without loss of generality,
 $|n_{j_1}(k)|{\geq}|n_{j_2}(k)|\geq 2^{1-B_{\max}}$), we  bound instead
 \begin{align}
 \nonumber
&\frac{|b_{\ell}(k{+}1)|}{\gamma(k)}{\stackrel{\triangle}{\leq}}\left\vert \frac{b_\ell(k)}{\gamma(k)}{-}d_{\ell}^+n_{\ell}(k)\right\vert
+\Bigr\vert n_{j_1}(k){+}n_{j_2}(k)\Bigr\vert \!\!
+\!\!\!\!\!\!\!\!\!\!\!\!\!\sum_{j\in \mathcal{N}_{\ell}^-,j\neq j_1,j_2}\!\!\!\!\!\!\!\!\!\!\!\!\vert n_{j}(k)\vert
\\
&
{\leq}\frac{|b_\ell(k)|}{\gamma(k)}-d_\ell^+|n_{\ell}(k)|
+\sum_{j\in \mathcal{N}_\ell^-}\vert n_{j}(k)\vert
-2|n_{j_2}(k)|.
\end{align}
In both cases,
since $|n_{j_2}(k)|{\geq}2^{1-B_{\max}}$, $|n_{j_1}(k)|{\geq}2^{1-B_{\max}}$ and
$\gamma(k)|n_{\ell}(k)|{\leq}|b_\ell(k)|/d_\ell^+$ (see \eqref{signal_wb_alg_ac}),
 we further bound
 \begin{align}
\label{bal_decr_2}
\!\!\!\!\!\frac{|b _\ell(k{+}1)|}{\gamma(k)}{\leq}\frac{|b_\ell(k)|}{\gamma(k)}{-}d_\ell^+|n_{\ell}(k)|
{+}\sum_{j\in \mathcal{N}_\ell^-}\vert n_{j}(k)\vert
{-}2^{2-B_{\max}}.
\end{align}
By summing
\eqref{qwer} (with strict equality if $\mathcal I\{\mathcal D_k\} = 0$)
 and \eqref{bal_decr_2} over $i\in \mathcal V$, it holds
\begin{align*}
\Vert {\bf b}(k+1)\Vert_1
\begin{cases}
= \Vert {\bf b}(k)\Vert_1, & \text{if }\mathcal I\{\mathcal D_k\} = 0 \\
\leq \Vert {\bf b}(k)\Vert_1-2\bar\gamma(k), & \text{otherwise,}
\end{cases}
\end{align*}
after noticing that $j\in\mathcal{N}_i^-\Leftrightarrow i\in\mathcal{N}_j^+$, hence
\begin{align*}
\sum_{i,j\in\mathcal{N}_i^-}\vert n_j(k) \vert
=
\sum_{j,i\in\mathcal{N}_j^+}\vert n_j(k)\vert
=\sum_{j=1}^N\vert n_j(k) \vert d_j^+.
\end{align*}
This completes the proof.\vspace{-0.3cm}

\subsection{Proof of Proposition \ref{lemma_Uk}}  \label{append:lemma_Uk}
\noindent\textbf{Property {(i)}:}  
 Note that $U(k){\geq}0$ since $\vert b_i(k) \vert{\geq}0$, for all $i \in \mathcal{V}$.
To show that it is upper bounded, we use
$U_h{\leq}U_1,\forall h{\geq}1$ and  
$\cup_h \mathcal{V}_h(k)\subseteq\mathcal V$, and write 
\begin{align*}
U(k)
&\leq \sum\limits_{h = 1}^{N-1}{U_n\sum\limits_{i \in \mathcal{V}_h(k)}{N}} 
\leq U_1N\sum\limits_{h = 1}^{N-1}{|\mathcal{V}_h(k)|} 
\leq N^{2(N-1)}.
\end{align*}
We prove that $U(k)$ is nondecreasing as by product of the proof of Property (ii), as given below.  

\noindent\textbf{Property {(ii)}:}  
Since $\bar\gamma(k+1)\leq \bar\gamma(k)$, we can lower bound
   $U(k+1)$ as
   \begin{align*}
   &U(k+1){\geq }
   \sum\limits_{h = 1}^{N-1}U_h
\sum\limits_{j \in \mathcal{V}_h(k{+}1)}\min\left\{\frac{\vert b_j(k+1) \vert}{\bar\gamma(k)},N\right\}.
   \end{align*}
     \textbf{Case 1: $\mathcal{I}(\mathcal{D}_k)=\mathcal{I}(\mathcal{U}_k)=0$.}  We have  $b_j(k+1)=b_j(k),\forall j$
      and $\mathcal{V}_h(k{+}1)\triangleq\mathcal{V}_h(k),\forall n$, 
   which implies 
   $U(k+1)\geq U(k)$. \textbf{Case 2: $\mathcal{I}(\mathcal{U}_k)=1$.} 
   From the discussion following \eqref{qwer},       \eqref{qwer} holds with equality:
   \begin{align}
   \label{vnm}
\!\!\!\!\!\frac{|b_i(k+1)|}{\gamma(k)}\stackrel{\triangle}{=} \frac{|b_i(k)|}{\gamma(k)}-d_i^+n_{i}(k)+\sum_{j\in \mathcal{N}_i^-}\vert n_{j}(k)\vert,\ \forall i.
\end{align}
    Moreover, $\exists i\in\mathcal V:n_i(k)\neq 0$; this
    implies that there exists a non-empty set of nodes that receive at least one update from their in-neighbors, defined as
   $$\mathcal R(k)=\{i\in\mathcal V: n_j(k)\neq 0,\exists j\in\mathcal N_i^-\}.$$
   It is straightforward to show that 
\begin{align}
\label{fghjk}
\mathcal R(k)\subseteq \mathcal V^+(k+1)\cup\mathcal V^-(k+1).
\end{align}
In fact, if $i\notin \mathcal V^+(k+1)\cup\mathcal V^-(k+1)$ (i.e., $b_i(k+1)=0$),
    it follows that $i\notin\mathcal V^+(k)\cup\mathcal V^-(k)$ (i.e., $b_i(k)=0$ and $n_{i}(k)=0$, cf. L.\ref{lemma_node_flow});
    therefore, setting $b_i(k+1)=b_i(k)=0$ and $n_i(k)=0$ in \eqref{vnm}, we find that 
    $n_{j}(k)=0,\ \forall j\in\mathcal{N}_i^-$,
    so that $i\notin\mathcal R(k)$ and \eqref{fghjk} follows.
    With this definition, let
    $$
    h^*=\min\{h\in\{1,2,\dots,N-1\}:|\mathcal V_h(k+1)\cap\mathcal R(k)|>1\}
    $$
    be the distance of the node closest to those of opposite sign of balance at $k+1$ to receive the update,
    and let $\ell\in\mathcal V_{h^*}(k+1)\cap\mathcal R(k)$ be one of  such nodes.
   Then, we have 
\begin{align}
\label{nik=0}
    n_i(k)=0,\ \forall i \in\mathcal V_h(k+1),\,\,\forall h\leq h^*.
\end{align}
In fact, if $n_i(k)\neq 0$ for some of such $i$, then $\exists j\in\mathcal N_i^+\cap\mathcal V_{h-1}(k+1)\cap\mathcal R(k)$ receiving the update, which contradicts the definition of $h^*$.
Reading \eqref{vnm} at $i{\in}\mathcal V_h(k{+}1), h{\leq}h^*$, yields 
\begin{align*}
\frac{|b_i(k+1)|}{\bar\gamma(k)}{=}\frac{|b_i(k)|}{\bar\gamma(k)}{+}\sum_{j\in \mathcal{N}_i^-}\frac{\vert n_{j}(k)\vert}{2^{1-B_{\max}}}{\geq }
 \frac{\vert b_i(k)\vert}{\bar\gamma(k)}{+}\mathcal I\{i{\in}\mathcal R(k)\}.
\end{align*}  
We can then further lower bound $U(k+1)$ as
       \begin{align}
          \label{nhdft2}
          \nonumber
   U(k{+}1){\geq }&
   \sum\limits_{h = 1}^{h^*}U_h\!\!\!\!\!\!\!
\sum\limits_{j \in \mathcal{V}_h(k{+}1)}\!\!\!\!\!\!\!\min\left\{\frac{\vert b_j(k) \vert}{\bar\gamma(k)},N\right\}\mathcal I\{(h,j){\neq}(h^*,\ell)\}
\\&
   + U_{h^*}\min\left\{ \frac{|b_i(k)|}{\bar\gamma(k)}+1,N\right\},
   \end{align}
   where we neglected the non-negative terms associated to 
   $\mathcal V_h(k+1),h>h^*$.
   To further bound this quantity, note that $n_\ell(k)=0$ (cf. \eqref{nik=0}),
   which, together with $B_i^{\rm (w)}(k)>0, \forall i$ for an update event to occur, implies $b_\ell(k)/\bar\gamma(k)\leq d_i^+\leq N-1$ (cf. \eqref{signal_wb_alg_ac}).
      Therefore $\min\left\{\vert b_\ell(k)\vert/\bar\gamma(k)+1,N\right\}=\min\left\{\vert b_\ell(k)\vert/\bar\gamma(k),N\right\}+1$.
     Finally, we use the fact that 
     $$
     0\geq 
     \sum\limits_{h = h^*+1}^{N-1}U_h
\sum\limits_{j\in\mathcal{V}_h(k+1)}\left[\min\left\{\frac{\vert b_j(k) \vert}{\bar\gamma(k)},N\right\}-N\right],
     $$
     yielding
     \begin{align}
     \nonumber
  \label{nhdft3}
   U(k+1){\geq }&
 \sum\limits_{j \in \mathcal{V}}
\min\left\{\frac{\vert b_j(k) \vert}{\bar\gamma(k)},N\right\}   \sum\limits_{h = 1}^{N-1}U_h\mathcal I\{j \in \mathcal{V}_h(k{+}1)\}\\&
  +U_{h^*}-N\sum\limits_{h = h^*+1}^{N-1}U_h|\mathcal{V}_h(k+1)|.
   \end{align}
Now, using $U_h{\leq}U_{h^*+1},\forall h{>}h^*$ and 
   $\cup_{h = h^*+1}^{N-1}\mathcal{V}_h(k{+}1)\subseteq\mathcal V\setminus\{\ell\}$, we obtain
\begin{align}
\label{boundU}
&   U_{h^*}-N\sum\limits_{h = h^*+1}^{N-1}U_h|\mathcal{V}_h(k+1)|
  \\&
   \geq
   U_{h^*}-NU_{h^*+1}\sum\limits_{h = h^*+1}^{N-1}|\mathcal{V}_h(k+1)|
   \geq 
   N^{2(N-h^*)-3}
   \geq
   N.
   \nonumber
\end{align}
   In the last inequality, we used the fact that $h^*{\leq}N-2$. In fact,
   if $h^*{=}N{-}1$, then \eqref{nik=0} implies that $n_i(k){=0},\forall i$, which contradicts the occurrence of the update event.   
   Finally, note that
   $\mathcal V^+(k){\subseteq}\mathcal V^+(k{+}1)$ and $\mathcal V^-(k){\subseteq} \mathcal V^-(k{+}1)$ (cf. L.\ref{lemma_node_flow}), hence
   $j{\in}\mathcal{V}_h(k)\Rightarrow j{\in}\cup_{m=1}^h\mathcal{V}_m(k{+}1)$, i.e.,
   $j$ gets closer to nodes of opposite sign.
   Together with $U_h{>}U_{h+1}$, it implies
   \begin{align}
   \label{nmq}
   \sum\limits_{h = 1}^{N-1}U_h\mathcal I\{j \in \mathcal{V}_h(k{+}1)\}
   \geq
   \sum\limits_{h = 1}^{N-1}U_h\mathcal I\{j \in \mathcal{V}_h(k)\}.
   \end{align}
   The desired result follows by using \eqref{boundU}-\eqref{nmq} in \eqref{nhdft3}.

\noindent\textbf{Property {(iii)}:} We prove it by contradiction. Assume that $t_1, t_2 \in [k,t_{l+1})$ such that $\mathcal U_{t_1}$ and $\mathcal U_{t_2}$ are two consecutive update events with $t_2 - t_1 > 2W-1$, and $\Vert {\bf b}(t)\Vert_1 \geq N^2\gamma(t),\forall t \in [k,t_{l+1})$. It follows that $\exists i \in \mathcal V$ such that
$$\vert b_i(t_1) \vert 
\geq N\gamma(t_1) 
\stackrel{(A.2)}{\geq} N\gamma(t)
> d_i^+\gamma(t), \forall t \in [t_1, t_1+2W-1],
$$
which implies that $\exists t \in [t_1, t_1+2W-1]:\mathcal I\{\mathcal U_t\}=1$ due to A.\ref{assump_Bw}, which contradicts the assumption that $\mathcal U_{t_1}$ and $\mathcal U_{t_2}$ are two consecutive update events since $t_1+2W-1 < t_2$. Hence, it follows from property {\it (ii)} that
\begin{align*}
U(k+2W-1) \geq U(k)+N, 
\end{align*} 
which proves property {\it (iii)}.
\hfill $\blacksquare$

\section{Auxiliary Results for Theorem \ref{thm_ac}} \label{append:ac}

\begin{lemma} \label{lemma_avg_presv}
Let $\{{\bf y}(k)\}_{k \in \mathbb{Z}_+}$ be the sequence generated by Algorithm \ref{alg_dac}, in the setting of T.\ref{thm_ac}. Then, it holds:
\begin{align*}
{\bf 1}^\top{\bf y}(k) = {\bf 1}^\top{\bf y}(0)
\end{align*}
\end{lemma}
\begin{IEEEproof}
From \eqref{update_y} and ${\bf L}^{+}(k)={\bf S}^+(k)-{\bf A}(k)$, it follows
\begin{align}
\mathbf 1^T{\bf y}(k+1){=}\mathbf 1^T{\bf y}(k){-}\alpha(k)\mathbf 1^T[{\bf S}^+(k){-}{\bf A}(k)]\mathbf x(k), \label{average_preserve}
\end{align}
so that the statement of the lemma readily follows after noticing that 
$\mathbf 1^T{\bf A}(k)={\bf S}^+(k)$.
\end{IEEEproof}

\begin{lemma}\label{lemma_y_global_bound}
Let $\{{\bf y}(k)\}_{k \in \mathbb{Z}_+}$ be the sequence generated by Algorithm \ref{alg_dac}, in the setting of T.\ref{thm_ac}. Then,
$$y_{i,\min}\leq
y_{i,\min}^{(k)}\leq
 y_i(k)
\leq  y_{i,\max}^{(k)}\leq y_{i,\max},\ \forall k{\in}\mathbb Z_+,
$$
where $y_{i,\max}{\triangleq}\lim_{t\to\infty}y_{i,\max}^{(t)}{<}\infty,y_{i,\min}{\triangleq}\lim_{t\to\infty}y_{i,\min}^{(t)}{>}\infty$, $q^*{\triangleq}\max\{|q_{\min}|,|q_{\max}|\}$, $S_{\max}{\triangleq}\sup_{k \in \mathbb{Z}_+}S_i^-(k){<}\infty$,
\begin{align}
\nonumber
&y_{i,\max}^{(k)}\triangleq
\max\{q_{\max},y_i(0)\}+\alpha(0) \Vert {\bf b}(0)\Vert_1 q^*\\
&+\alpha(0)S_{\max}\left(q_{\max}-q_{\min}\right)+q^*\sum_{t=0}^{k-1}{\alpha(t)\vert b_i(t)\vert},
\\
\nonumber
&
y_{i,\min}^{(k)}\triangleq
\min\{q_{\min},y_i(0)\}-\alpha(0) \Vert {\bf b}(0)\Vert_1 q^*\\
&-\alpha(0)S_{\max}\left(q_{\max}-q_{\min}\right)-q^*\sum_{t=0}^{k-1}{\alpha(t)\vert b_i(t)\vert}.
\end{align}
\end{lemma}
\begin{IEEEproof}
We first show that $|y_{i,\max}|,|y_{i,\min}|{<}\infty$. From T.\ref{thm_wb}, we know that ${\bf A}(k)$ is bounded for all $k{\in}\mathbb{Z}_+$, which implies $S_i^-(k)\leq S_{\max} < \infty$. On the other hand, since $\vert b_i(t)\vert{\leq}\Vert {\bf b}(t)\Vert_1{=}\mathcal{O}\left(\frac{1}{t}\right)$ and $\sum_{t \in \mathbb{Z}_+}{\alpha(t)^2}{<}\infty$, one can verify using Cauchy-Schwarz inequality that $\sum_{t=0}^\infty{\alpha(t)\vert b_i(t)\vert}{<} \infty$ and thus $|y_{i,\max}|,|y_{i,\min}|{<}\infty$.
By inspection, it is also clear that 
$y_{i,\max}^{(k)}\leq y_{i,\max}$ and $y_{i,\min}^{(k)}\geq y_{i,\min},\forall k$.
We now prove $y_i(k){\in}[y_{i,\min}^{(k)},y_{i,\max}^{(k)}]$ by induction.
Clearly, it for $k=0$. Now, assume it holds for some $k\geq 0$, we prove that this implies $y_i(k{+}1){\in}[y_{i,\min}^{(k{+}1)},y_{i,\max}^{(k{+}1)}]$
 (induction step). We have:
\begin{enumerate}
\item If $y_i(k)\leq  q_{\max}$, then
\begin{align*}
&y_i(k+1) = y_i(k)+\alpha(k)b_i(k)x_i(k)\\
&+\alpha(k)\sum\limits_{j \in \mathcal{N}_i^-}{a_{ij}(k)\Big(x_j(k)-x_i(k)\Big)} \\
&\leq q_{\max}+\alpha(k) \vert b_i(k)\vert q^* + \alpha(k)S_i^-(k) \left(q_{\max}-q_{\min}\right) \\
&{\leq} \max\{q_{\max},y_i(0)\}{+}\alpha(0) \Vert {\bf b}(0)\Vert_1 q^*\\&{+}\alpha(0)S_{\max}\left(q_{\max}{-}q_{\min}\right)
=y_{i,\max}^{(0)}\leq y_{i,\max}^{(k+1)}.
\end{align*}
\item If $y_i(k) > q_{\max}$, then $x_i(k) = q_{\max}$, so that \eqref{update_consensus} yields
\begin{align*}
y_i(k+1)&\leq y_i(k)+\alpha(k)b_i(k)q_{\max} \\
&\leq y_{i,\max}^{(k)}+\alpha(k)|b_i(k)|q^*
=y_{i,\max}^{(k+1)}.
\end{align*}
\item
Similarly, if $y_i(k)\geq q_{\min}$ then
\begin{align*}
&y_i(k{+}1){\geq}\min\{q_{\min},y_i(0)\}{-}\alpha(0) \Vert {\bf b}(0)\Vert_1 q^*
\\&{-}\alpha(0)S_{\max}\left(q_{\max}{-}q_{\min}\right)
=y_{i,\min}^{(0)}\geq y_{i,\min}^{(k+1)}.
\end{align*}
\item If $y_i(k) < q_{\min}$, then $x_i(k) = q_{\min}$ so that
\eqref{update_consensus} yields
\begin{align*}
y_i(k+1) &\geq y_i(k)-\alpha(k) \vert b_i(k)\vert q^*
\\&
\geq y_{i,\min}^{(k)}-\alpha(k) \vert b_i(k)\vert q^*=y_{i,\min}^{(k+1)}.
\end{align*}
\end{enumerate}\vspace{-0.4cm}
\end{IEEEproof}

\begin{lemma} \label{lemma_ac_dec_term}
Let $\{{\bf y}(k)\}_{k \in \mathbb{Z}_+}$ be the sequence generated by Algorithm \ref{alg_dac}, in the setting of T.\ref{thm_ac}. Then,
\begin{align}\label{eq:lower_bound_negative_term}
{\bf y}(k)^\top{\bf L}^+(k)\tilde{\bf y}(k) \geq \xi_1 V({\bf y}(k)) - \xi_2 \Vert {\bf b}(k)\Vert_1,
\end{align}
for some finite constants $\xi_1, \xi_2 > 0$.
\end{lemma}
\begin{IEEEproof}
Let $\hat{\bf e}(k){=}{\bf y}(k){-}\tilde{\bf y}(k)$ be the saturation error,
\\\noindent ${\bf S}^{\pm}\!(k){=}{\rm diag}\!\left\{\!S_i^{\pm}\!(k),\forall i\right\}$, 
${\bf B}(k){=}{\rm diag}\!\left\{\!\mathbf b(k)\right\}{=}{\bf S}^-\!(k){-}{\bf S}^+\!(k)$,
\\\noindent
${\bf L}^{\pm}\!(k){=}{\bf S}^{\pm}\!(k){-}{\bf A}\!(k)$, 
 ${\bf L}(k){=}[{\bf S}^+\!(k){+}{\bf S}^-\!(k)]{-}[{\bf A}\!(k){+}{\bf A}\!(k)^\top]$. The proof contains three steps:
\begin{itemize}
\item {\bf Step 1:} We will lower bound ${\bf y}(k)^\top{\bf L}^+(k)\tilde{\bf y}(k)$ as
\begin{align}
&{\bf y}(k)^\top{\bf L}^+(k)\tilde{\bf y}(k) 
\geq -{\bf y}(k)^\top{\bf B}(k)\tilde{\bf y}(k)  \nonumber\\
&+\frac{1}{2}\tilde{\bf y}(k)^\top{\bf B}(k)\tilde{\bf y}(k)+\frac{1}{2}\tilde{\bf y}(k)^\top{\bf L}(k)\tilde{\bf y}(k).\label{y_2nd_term_eq1}
\end{align}

\item {\bf Step 2:} we will show that the last term of the RHS in Step 1 satisfies, for some $\xi_4 > 0$,
\begin{align}
\tilde{\bf y}(k)^\top{\bf L}(k)\tilde{\bf y}(k) \geq \xi_4 V({\bf y}(k)).\label{2nd_term_lb}
\end{align}
\item {\bf Step 3:} by combining the above results, we will show that, for some constants $\xi_1,\xi_2 > 0$,
\begin{align*}
{\bf y}(k)^\top{\bf L}^+(k)\tilde{\bf y}(k) \geq \xi_1 V({\bf y}(k)) - \xi_2 \Vert {\bf b}(k)\Vert_1.
\end{align*}
\end{itemize}
In the following, we provide detailed derivations of each step.
\begin{itemize}
\item {\bf Step 1:} 
It is easy to show that
\begin{align*}
\!\!\!\!\!\!{\bf y}(k)^\top{\bf L}^+(k)\tilde{\bf y}(k) {=} 
 {-}{\bf y}(k)^\top{\bf B}(k)\tilde{\bf y}(k){+}{\bf y}(k)^\top{\bf L}^-(k)\tilde{\bf y}(k).
\end{align*}
The term ${\bf y}^\top(k){\bf L}^-(k)\tilde{\bf y}(k)$ can be lower bounded as
\begin{align*}
&\!\!{\bf y}(k)^\top{\bf L}^-(k)\tilde{\bf y}(k)
{=} \hat{\bf e}(k)^\top{\bf L}^-(k)\tilde{\bf y}(k){+}\tilde{\bf y}(k)^\top{\bf L}^-(k)\tilde{\bf y}(k) \\
&\!\!\overset{(a)}{\geq} \tilde{\bf y}(k)^\top{\bf L}^-(k)\tilde{\bf y}(k)
{=} \frac{1}{2}\tilde{\bf y}(k)^\top\left[{\bf L}^-(k){+}{\bf L}^-(k)^\top\right]\tilde{\bf y}(k) \\
&\!\!= \frac{1}{2}\tilde{\bf y}(k)^\top{\bf B}(k)\tilde{\bf y}(k)+\frac{1}{2}\tilde{\bf y}(k)^\top{\bf L}(k)\tilde{\bf y}(k),
\end{align*}
where $(a)$ comes from the fact that
\begin{align*}
&\hat{\bf e}(k)^\top{\bf L}^-(k)\tilde{\bf y}(k)=\hat{\bf e}(k)^\top[{\bf S}^-(k)-{\bf A}(k)]\tilde{\bf y}(k) \\
&= \sum_{i=1}^N{\Big[\hat e_i(k)\sum_{j=1}^N{a_{ij}(k)\left( \tilde y_i(k)-\tilde y_j(k)\right)} \Big]} \geq 0,
\end{align*}
where the last inequality comes from the fact that (i) if $y_i(k) \in [q_{\min}, q_{\max}]$, then $\hat e_i(k)=0$; (ii) if $y_i(k) > q_{\max}$, then $\hat e_i(k) > 0$ and $\tilde y_i(k)-\tilde y_j(k) = q_{\max}-\tilde y_j(k) \geq 0, \forall j \in \mathcal V$; and (iii) if $y_i(k) < q_{\min}$, then $\hat e_i(k) < 0$ and $\tilde y_i(k)-\tilde y_j(k) = q_{\min}-\tilde y_j(k) \leq 0, \forall j \in \mathcal V$.

\item {\bf Step 2:} 
First, one can verify that
\begin{align*}
&\tilde{\bf y}(k)^\top{\bf L}(k)\tilde{\bf y}(k) {=} \frac{1}{2}\sum\limits_{i,j=1}^N{[a_{ij}(k){+}a_{ji}(k)]\left[\tilde{y}_i(k){-}\tilde{y}_j(k)\right]^2}.
\end{align*}
Let $i^*{\in}\arg\max_i \{y_i(k)\}$, $j^*{\in}\arg\min_i \{y_i(k)\}$, $j^*\neq i^*$.
Note that ${ y}_{i^*}(k) \geq \bar y(0) \geq { y}_{j^*}(k)$ to preserve the average (L.\ref{lemma_avg_presv}).
 Since $\mathcal{G}$ is strongly connected, there exists a path from $i^*$ to $j^*$. Let $\{i_1,\cdots, i_p\}$ be the set of nodes in the shortest path from $i^*$ to $j^*$, with $i_1{=}i^*, i_p{=}j^*$ and $i_{n+1}{\in}\mathcal{N}_{i_n}^+, \forall n \in [1,p-1]$.
We have 
\begin{align}
&\tilde{\bf y}(k)^\top{\bf L}(k)\tilde{\bf y}(k)
{=}\frac{1}{2}\sum\limits_{i,j=1}^N{[a_{ij}(k){+}a_{ji}(k)]\left[\tilde{y}_i(k){-}\tilde{y}_j(k)\right]^2} \nonumber \\
&\geq \frac{1}{2}\sum\limits_{l=1}^{p-1}{[a_{i_li_{l+1}}(k)+a_{i_{l+1}i_l}(k)]\left[\tilde{y}_{i_l}(k)-\tilde{y}_{i_{l+1}(k)}\right]^2} \nonumber \\
&\overset{(a)}\geq \frac{a_{\min}}{2} \sum\limits_{l=1}^{p-1}{\left[\tilde y_{i_l}(k)-\tilde y_{i_{l+1}(k)}\right]^2} \nonumber \\
&\overset{(b)}\geq \frac{a_{\min}}{2(p-1)}\Big\{\sum\limits_{l=1}^{p-1}\left[{\tilde y_{i_l}(k)-\tilde y_{i_{l+1}}(k)}\right]\Big\}^2 \nonumber\\
&\geq\frac{a_{\min}}{2(N-1)}\left[\tilde y_{i^*}(k)-\tilde y_{j^*}(k)\right]^2, \label{lb_ytLyt}
\end{align}
where $(a)$ follows from $a_{i_{l+1}i_l}(k)\geq a_{\min},\forall l \in [1,p), \forall k \in \mathbb Z_+$ (T.\ref{thm_wb}$(iii)$);
 $(b)$ comes from Cauchy-Schwarz inequality.
To further bound this quantity, note that
 \begin{align}
 \nonumber
 &\frac{1}{N}V(\mathbf y)=\frac{1}{N}\sum_i[{y}_i(k)-\bar{y}(0)]^2
 \\&
\leq \max_{i\in\{i^*,j^*\}}[{y}_{i}(k)-\bar{y}(0)]^2
\leq [y_{i^*}(k)-y_{j^*}(k)]^2.
\label{bfntr}
\end{align}
On the other hand, since the consensus algorithm preserves the average,
it follows
\begin{align}
y_{i^*}(k)-\bar{y}(0)\leq (N-1)(\bar{y}(0)-y_{j^*}(k))
\label{L11_ineq2}
\\
\bar{y}(0)-y_{j^*}(k)\leq (N-1)(y_{i^*}(k)-\bar{y}(0))
,
\end{align}
so that the first inequality in \eqref{bfntr} is upper bounded as
 $$\frac{1}{N}V(\mathbf y)
\leq (N-1)^2\min_{i\in\{i^*,j^*\}}[{ y}_{i}(k)-\bar{y}(0)]^2.
$$
Consider the following two cases:
\begin{enumerate}[(i)]
\item $y_i(k){\in}[q_{\min},q_{\max}], \forall i$, so that $\tilde{y}_i(k){=}y_i(k), \forall i$ and 
\begin{align*}
&\tilde y_{i^*}(k)- \tilde y_{j^*}(k)
= y_{i^*}(k)-y_{j^*}(k)\geq \frac{1}{\sqrt{N}}\sqrt{V(\mathbf y)}.
\end{align*}
\item $y_{i^*}(k){>}q_{\max}$ ($y_{j^*}(k){<}q_{\min}$ can be solved similarly) so that $\tilde y_{i^*}(k){=}q_{\max}$: 
since $\tilde y_{j^*}(k){\leq}\max\{y_{j^*}(k),q_{\min}\}$, using \eqref{L11_ineq2}
and $\bar y(0)\leq q_{\max}$ it follows
\begin{align}
\nonumber
&\tilde y_{i^*}(k)-\tilde y_{j^*}(k){\geq }
\min\Big\{\frac{y_{i^*}(k)-\bar{y}(0)}{N-1},q_{\max}-q_{\min}\Big\}
\\&
\geq 
\min\Big\{\frac{\sqrt{V(\mathbf y)}}{\sqrt{N}(N-1)^2},q_{\max}-q_{\min}\Big\}.
\end{align}
\end{enumerate}

From (i), (ii) and \eqref{lb_ytLyt}, there exists some $\xi_4 > 0$ such that
\begin{align*}
&\tilde{\bf y}(k)^\top{\bf L}(k)\tilde{\bf y}(k) \geq \frac{a_{\min}}{2(N-1)}\left(\tilde y_{i^*}(k)-\tilde y_{j^*}(k)\right)^2 \\
&\!\!\!\!\!\!\geq \frac{a_{\min}/2}{(N{-}1)}\min\left\{\frac{V\left({\bf y}(k)\right)}{N(N{-}1)^4},(q_{\max}{-}q_{\min})^2\right\}{\geq}\xi_4 V({\bf y}(k)),
\end{align*}
since $y_i(k)$ and thus $V({\bf y}(k))$ is bounded (L.\ref{lemma_y_global_bound}).

\item {\bf Step 3:} Let $y^*= \max_i \{\max\{\vert y_{i,\max} \vert, \vert y_{i,\min} \vert\}\}$. By combining \eqref{y_2nd_term_eq1} and \eqref{2nd_term_lb}, we get
\begin{align*}
&\!\!\!\!{\bf y}(k)^\top{\bf L}^+(k)\tilde{\bf y}(k)\\
&\!\!\!\!\!\!\geq -{\bf y}(k)^\top{\bf B}(k)\tilde{\bf y}(k)+\frac{1}{2}\tilde{\bf y}(k)^\top{\bf B}(k)\tilde{\bf y}(k)+\xi_1  V({\bf y}(k))\\
&\!\!\!\!\!\!\geq \!{-}\!\sum_{i=1}^N{\vert b_i(k)y_i(k)\tilde{y}_i(k)\vert}{-}\frac{1}{2}\!\sum_{i=1}^N{\vert b_i(k)\tilde{y}_i(k)^2 \vert}{+}\xi_1 V({\bf y}(k))\\
&\!\!\!\!\!\!\overset{(a)}\geq 
-q^*\left(y^*+\frac{1}{2}q^*\right) \Vert {\bf b}(k)\Vert_1+\xi_1 V({\bf y}(k)),
\end{align*}
with $\xi_1{=}\xi_4/2{>}0$, where $(a)$ comes from the facts $\vert y_i(k) \vert \leq y^*$, $\vert \tilde{y}_i(k) \vert \leq q^*$, and $\Vert {\bf b}(k)\Vert_1=\sum_{i=1}^N{\vert b_i(k)\vert}$.
\end{itemize}
\end{IEEEproof}
\vspace{-0.2cm}


\bibliographystyle{IEEEtran}
\bibliography{IEEEabrv,ref}	\vspace{-0.4cm}
\begin{IEEEbiography}
   [{\includegraphics[width=1in,height=1.25in,clip,keepaspectratio]{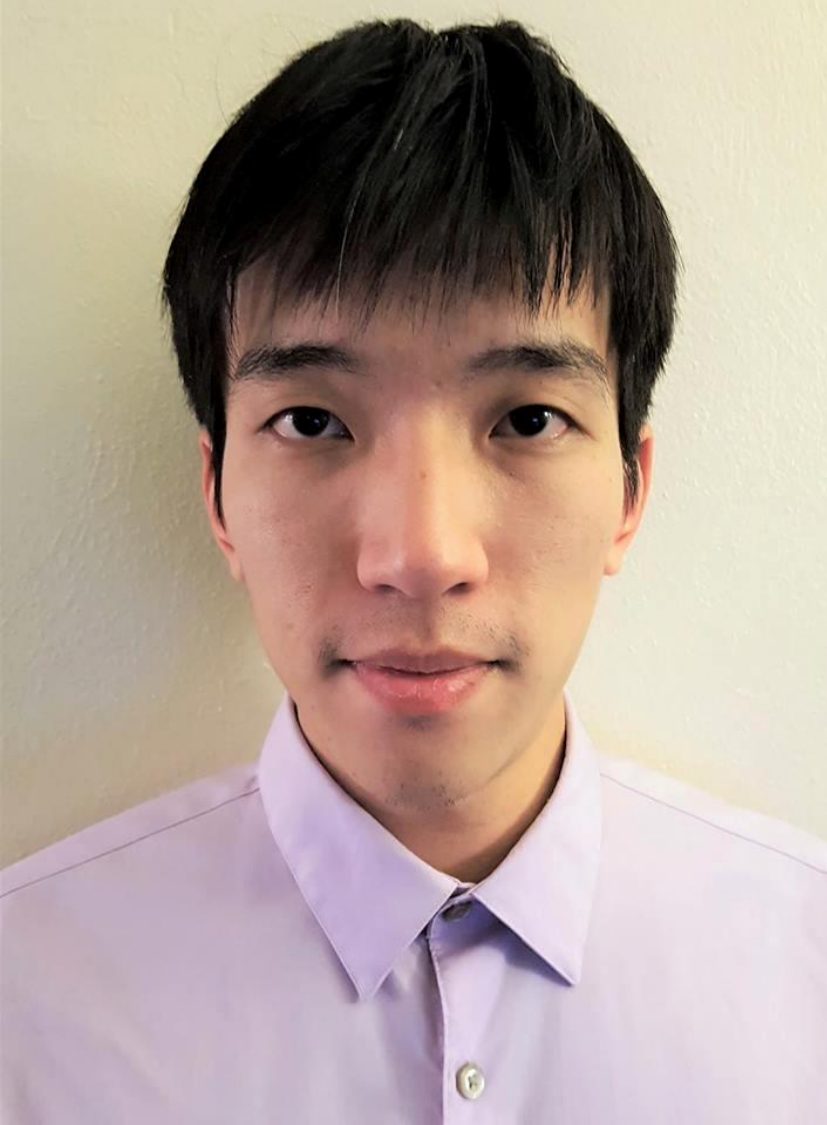}}]
    {Chang-Shen Lee} received the B.Sc. degree in the electrical engineering and computer science honor program and the M.Sc. in communications engineering from National Chiao Tung University, Hsinchu, Taiwan, in 2012 and 2014, respectively. He joined the Department of Electrical and Computer Engineering, Purdue University, West Lafayette, IN, in 2016, where he is currently pursuing the Ph.D. degree. From 2015 to 2016, he was a Research Assistant with the Research Center for Information Technology Innovation, Academia Sinica, Taiwan. His research interests include distributed computing and optimization.\vspace{-0.5cm}
\end{IEEEbiography}

\begin{IEEEbiography}
    [{\includegraphics[width=1in,height=1.25in,clip,keepaspectratio]{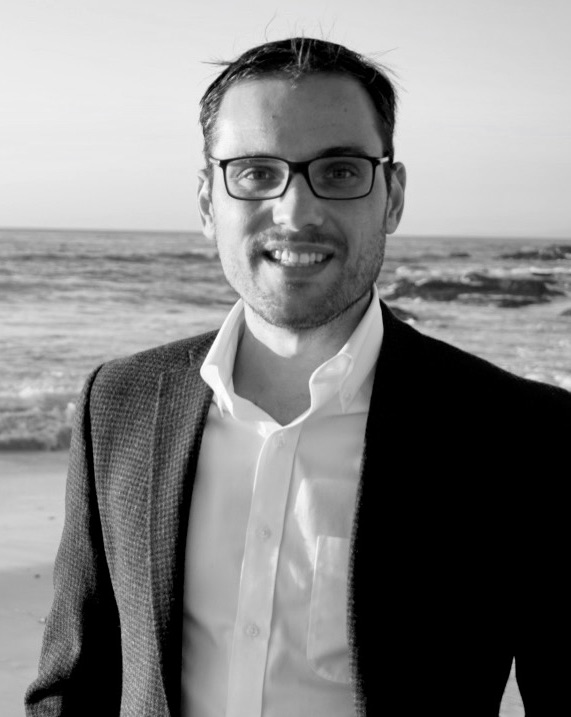}}]{Nicolo Michelusi}(S'09, M'13, SM'18) received the B.Sc. (with honors), M.Sc. (with honors) and Ph.D. degrees from the University of Padova, Italy, in 2006, 2009 and 2013, respectively, and the M.Sc. degree in Telecommunications Engineering from the Technical University of Denmark in 2009, as part of the T.I.M.E. double degree program. He was a post-doctoral research fellow at the Ming-Hsieh Department of Electrical Engineering, University of Southern California, USA, in 2013-2015. He is currently an Assistant Professor at the School of Electrical and Computer Engineering at Purdue University, IN, USA. His research interests lie in the areas of 5G wireless networks, millimeter-wave communications, stochastic optimization, distributed optimization. Dr. Michelusi serves as Associate Editor for the IEEE Transactions on Wireless Communications, and as a reviewer for several IEEE Transactions.\vspace{-0.5cm}
\end{IEEEbiography}
\begin{IEEEbiography}[{\includegraphics[width=1in,height=1.25in,clip,keepaspectratio]{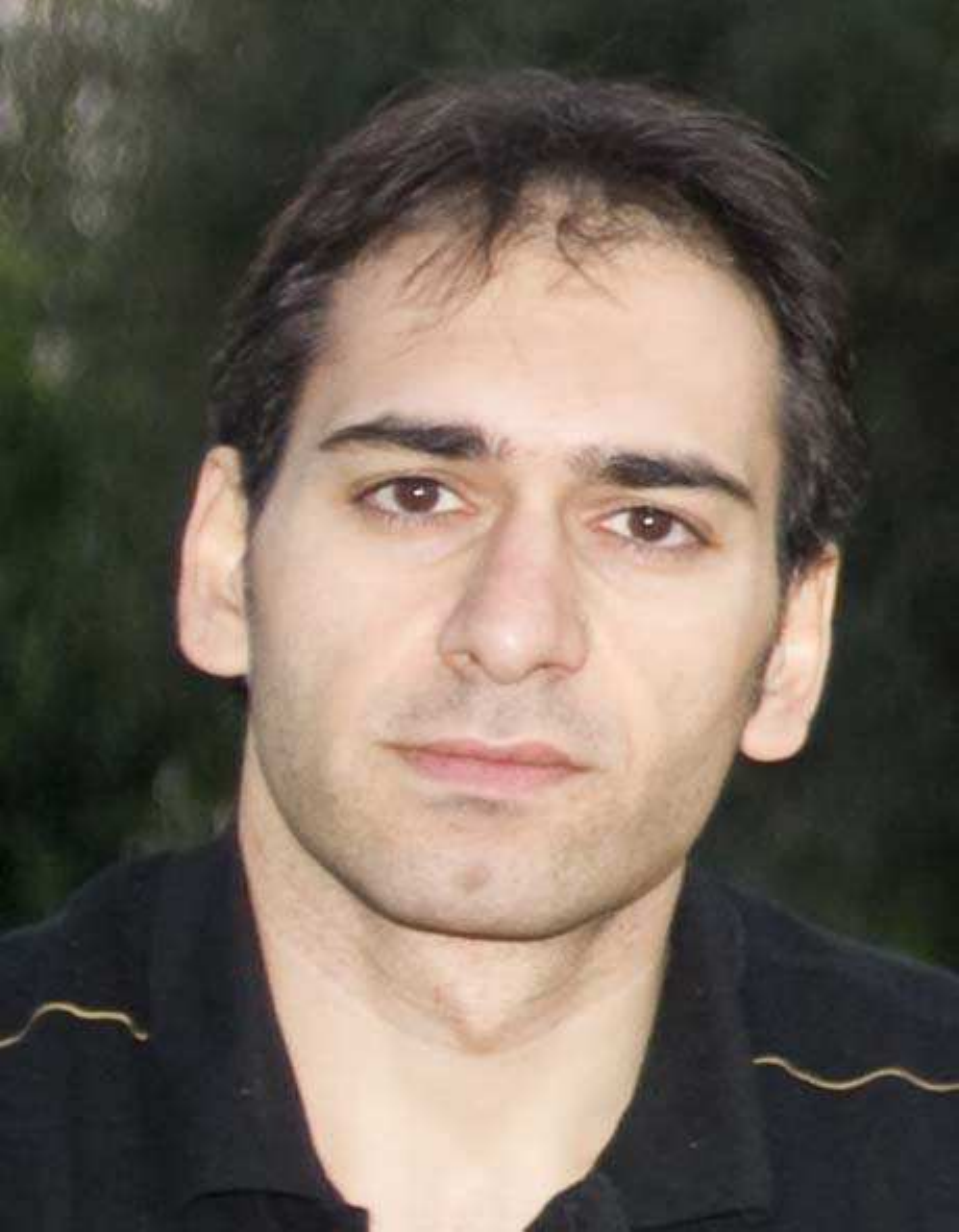}}]{Gesualdo Scutari}
(S'05-M'06-SM'11) received the Electrical Engineering and Ph.D. degrees (both with honors) 
from the University of Rome ``La Sapienza,'' Rome, Italy, in 2001 and 2005, respectively. 
He is the Thomas and Jane Schmidt Rising Star Associate Professor  with the School of Industrial Engineering, Purdue University, West Lafayette, IN, USA.
He had previously held several research appointments, namely, at the University of California at 
Berkeley, Berkeley, CA, USA; Hong Kong University of Science and Technology, Hong Kong; 
and University of Illinois at Urbana-Champaign, 
Urbana, IL, USA.  His research interests include continuous and distributed optimization,   equilibrium programming, and their applications to signal processing and 
machine learning. 
He is a Senior Area Editor of the IEEE Transactions On Signal Processing and an Associate Editor of Siam J. on Optimization; he served as an Associate Editor of the IEEE Transactions on Signal and Information Processing over Networks; the IEEE Transactions on Signal Processing, and the IEEE Signal Processing Letters. 

 He served on the IEEE Signal Processing Society Technical Committee on Signal Processing for Communications (SPCOM). 
He was the recipient of the 2006 Best Student Paper Award at the IEEE ICASSP 2006, the 2013 NSF  CAREER Award,  the 2015 Anna Maria Molteni Award for Mathematics and Physics, and the 2015 IEEE Signal Processing Society Young Author Best Paper Award.\vspace{-0.3cm}
\end{IEEEbiography}

\end{document}